\documentclass[sn-mathphys-num]{sn-jnl}

\usepackage{microtype}

\usepackage{graphicx}%
\usepackage{multirow}%
\usepackage{amsmath,amssymb,amsfonts}%
\usepackage{amsthm}%
\usepackage{mathrsfs}%
\usepackage[title]{appendix}%
\usepackage{xcolor}%
\usepackage{textcomp}%
\usepackage{manyfoot}%
\usepackage{booktabs}%
\usepackage{algorithm}%
\usepackage{algorithmicx}%
\usepackage{algpseudocode}%
\usepackage{listings}
\usepackage{mathtools}   
\usepackage{physics}
\usepackage{stmaryrd}    
\usepackage{tikz-cd}     
\usepackage{esint}       

\usepackage{subcaption}  
\usepackage{orcidlink}
\usepackage[most]{tcolorbox}
\usepackage{aliascnt}
\usepackage{dsfont} 
\usepackage{tabularray}
\usepackage{geometry}
\usepackage{bm}
\definecolor{Gray}{gray}{0.85}
\definecolor{LightCyan}{rgb}{0.8314, 0.9373, 0.9882}
\definecolor{BlueImperial}{rgb}{0,0.243,0.455}
\definecolor{Green}{rgb}{0.1333, 0.5451, 0.1333}
\definecolor{burgundy}{rgb}{0.5, 0.0, 0.13}
\definecolor{burntsienna}{rgb}{0.91, 0.45, 0.32}
\definecolor{chestnut}{rgb}{0.8, 0.36, 0.36}
\definecolor{cinnamon}{rgb}{0.82, 0.41, 0.12}
\definecolor{darkterracotta}{rgb}{0.8, 0.31, 0.36}
\definecolor{brightmaroon}{rgb}{0.76, 0.13, 0.28}
\definecolor{burlywood}{rgb}{0.87, 0.72, 0.53}
\definecolor{cinereous}{rgb}{0.6, 0.51, 0.48}
\definecolor{beige}{rgb}{0.972549, 0.901961, 0.768627}
\definecolor{champagne}{rgb}{0.97, 0.91, 0.81}
\definecolor{deepchampagne}{rgb}{0.98, 0.84, 0.65}
\definecolor{lightsalmonpink}{rgb}{1.0, 0.6, 0.6}
\definecolor{amber}{rgb}{1.0, 0.49, 0.0}
\definecolor{brightpink}{rgb}{1.0, 0.0, 0.5}
\definecolor{coquelicot}{rgb}{1.0, 0.22, 0.0}
\definecolor{antiquewhite}{rgb}{0.98, 0.92, 0.84}
\definecolor{cosmiclatte}{rgb}{1.0, 0.97, 0.91}

\newcommand\primarycolor{burgundy}     
\newcommand\secondarycolor{cinereous}  
\newcommand\boxcolorprim{lightsalmonpink} 
\newcommand\boxcolorsec{beige} 
\newcommand\urlcolor{blue} 
\newcommand\citecolor{amber} 

\hypersetup{%
        colorlinks,
        breaklinks=true,
        plainpages=false,%
        citecolor=\citecolor,
        linkcolor=\primarycolor,
        urlcolor=\urlcolor,
        bookmarksopen=true,%
        bookmarksnumbered=false,%
        bookmarksdepth=5%
}

\newtheoremstyle{thmstyleone}
{18pt plus2pt minus1pt}
{18pt plus2pt minus1pt}
{\itshape}
{0pt}
{\color{\primarycolor}\bfseries}
{}
{.5em}
{}

\numberwithin{equation}{section}

\theoremstyle{thmstyleone}
\newtheorem{theorem}{Theorem}[section]

\newaliascnt{corollary}{theorem}
\newtheorem{corollary}[corollary]{Corollary}
\aliascntresetthe{corollary}

\newaliascnt{lemma}{theorem}
\newtheorem{lemma}[lemma]{Lemma}
\aliascntresetthe{lemma}

\newaliascnt{definition}{theorem}
\newtheorem{definition}[definition]{Definition}
\aliascntresetthe{definition}

\newaliascnt{proposition}{theorem}
\newtheorem{proposition}[proposition]{Proposition}
\aliascntresetthe{proposition}

\newaliascnt{remark}{theorem}
\newtheorem{remark}[remark]{Remark}
\aliascntresetthe{remark}

\newaliascnt{notation}{theorem}
\newtheorem{notation}[notation]{Notation}
\aliascntresetthe{notation}

\newaliascnt{assumptions}{theorem}
\newtheorem{assumptions}[assumptions]{Assumptions}
\aliascntresetthe{assumptions}

\newaliascnt{example}{theorem}
\newtheorem{example}[example]{Example}
\aliascntresetthe{example}

\newaliascnt{conjecture}{theorem}
\newtheorem{conjecture}[conjecture]{Conjecture}
\aliascntresetthe{conjecture}

\newcommand{\Eq}[2]{\begin{equation}\label{#1}\begin{aligned}#2 \end{aligned}\end{equation}}

\newcommand{\theo}[2]{\rbox{\begin{theorem}\label{#1} #2 \end{theorem}}}
\newcommand{\coro}[2]{\rbox{\begin{corollary}\label{#1} #2 \end{corollary}}}
\newcommand{\lem}[2]{\bbox{\begin{lemma}\label{#1} #2 \end{lemma}}}
\newcommand{\defi}[2]{\bbox{\begin{definition}\label{#1} #2 \end{definition}}}
\newcommand{\prop}[2]{\bbox{\begin{proposition}\label{#1} #2 \end{proposition}}}
\newcommand{\rem}[2]{\begin{remark}\label{#1} #2 \end{remark}}
\newcommand{\nota}[2]{\bbox{\begin{notation}\label{#1} #2 \end{notation}}}

\newcommand{\conj}[2]{\bbox{\begin{conjecture}\label{#1} #2 \end{conjecture}}}

 \newcommand{\bbox}[1]{\begin{tcolorbox}[arc=0mm,oversize,left=15pt, right=17pt,colback=\boxcolorsec!20!white,colframe=\secondarycolor!100!white]#1\end{tcolorbox}}
\newcommand{\rbox}[1]{\begin{tcolorbox}[ arc=0mm,oversize,left=15pt, right=17pt,colback=\boxcolorprim!10!white,colframe=\primarycolor!100!white]#1\end{tcolorbox}} 

\renewcommand{\geq}{\geqslant}
\renewcommand{\leq}{\leqslant}

\title{\textbf{Critical phenomena in 2D Coulomb gas: Birth of a gap}}
\author*[1,2]{\fnm{Matthias} \sur{Allard} \orcidlink{0000-0002-5682-424X}}\email{m.allard@unimelb.edu.au}
\author*[1,2]{\fnm{Sampad} \sur{Lahiry}}\email{sampad.lahiry@kuleuven.be}

\affil[1]{\textit{\orgdiv{School of Mathematics and Statistics}, \orgname{University of Melbourne}, \orgaddress{\street{813 Swanston Street}, \city{Parkville, Melbourne}, \postcode{3010}, \state{Victoria}, \country{Australia}}}}
\affil[2]{\textit{\orgdiv{Department of Mathematics}, \orgname{Katholieke Universiteit Leuven}, \orgaddress{\street{Celestijnenlaan 200 B bus 2400}, \city{3001 Leuven},  \country{Belgium}}}}

\begin{document}

\abstract{
We investigate a  family of radially symmetric Coulomb gas systems at inverse temperature $\beta = 2$. The family is characterised by the property that the density of the equilibrium measure vanishes on a ring at radius $r_*$, which lies strictly inside the droplet. The large-$n$ expansion of the logarithm of the partition function is obtained up to the constant term and features a novel term of order $n^{1/4}$. We perform a double scaling limit of the correlation kernel at the $n^{1/4}$ scale and obtain a new limiting kernel in the bulk, which differs from the well-known Ginibre kernel.
}



\maketitle
\tableofcontents

\section{Introduction}
The study of two dimensional Coulomb gas is the study of a point process in the complex plane for which the probability measure is given by
\begin{equation}\label{eq:2D CG}
   d\mu_n^{(\beta)}(z)=\frac{1}{Z_n^{(\beta)}}\prod_{1\leq j<k\leq n}|z_j-z_k|^{\beta}\prod_{j=1}^{n}e^{-\frac{\beta n}{2}Q(z_j)}dA(z_j), 
\end{equation}
where $dA(z):=dxdy$, $(z=x+iy)$, is the flat Lebesgue measure on $\mathbb{C}$, $Z_n^{(\beta)}$ is a normalisation constant and $\beta>0$ is a parameter known in the physics literature as the inverse temperature. 
This models $n$ particles in interaction in a 2D plane exhibiting some logarithmic pairwise repulsion. This is best seen when rewriting the measure in the following form
\Eq{}{
     d\mu_n^{(\beta)}(z)=\frac{1}{Z_n^{(\beta)}}e^{-\frac{\beta }{2}H_n(z)}\prod_{j=1}^{n}dA(z_j),\quad H_{n}(z):=\sum_{1\leq j<k\leq n}\log( \frac{1}{|z_j-z_k|})+n\sum_{j=1}^{n}Q(z_j).
}
The particles are plunged in an external potential, given by the function 
\Eq{}{Q: \mathbb C \rightarrow \mathbb{R}\cup \{+\infty\}.}
The potential $Q$ is generally assumed to have sufficient growth at infinity to counterbalance the repulsion of the particles and confine a macroscopic number---proportional to $n$---of them in a bounded region of the complex plane. 

The point like particles $z=(z_1,\ldots,z_n)$ can be interpreted as charged particles at any given inverse temperature $\beta>0$. For the particular value of $\beta=2$, they can also be seen as the eigenvalues of a random normal matrix $M$ whose distribution is given by
\begin{equation}\label{normalmodel}
   \frac{1}{\mathcal{Z}_n^{(2)}}e^{-n \Tr(Q(M))}dM,
\end{equation}
with $dM$ the induced measure on the set of normal matrices and $\mathcal{Z}_n^{(2)}$ a different normalisation constant.

A wide variety of variations arise depending not only on the choice of the potential or family of potentials, but also on the value of $\beta$, e.g. \cite{BKS23,LS18,JV23,GKL24,ACC23,DGZ14,C24,C23,B25,AFLS25,CFTW15,Byun2024,N25}.
To make the above idea precise, one finds that as $n\rightarrow\infty$ the points distributed according to \eqref{eq:2D CG} accumulate in a 2-dimensional domain known as the droplet $\mathbb S$ according to its equilibrium measure $\mu_Q$. The measure $\mu_Q$ is characterised by the weighted logarithmic energy problem \cite{ST97}. In particular, for a  probability measure, define the weighted logarithmic energy by 
\begin{equation}\label{eq: energy}
    I_{Q}[\mu]:=\int_{\mathbb C^{2}}\log( \frac{1}{|z-w|})d\mu(z)d\mu(w)+\int_{\mathbb C} Q(w)d\mu(w).
\end{equation}
Then $\mu_Q$ minimises \eqref{eq: energy} over all probability measures supported on the complex plane. From Frostman's Theorem \cite{ST97}, it is well known that 
\begin{equation}\label{eq:eqmeas}
    d\mu_Q=\frac{1}{4\pi}\Delta Q\, \mathbf{1}_{\mathbb S}\, dA,
\end{equation}
where $\Delta:=4\partial \overline{\partial}$ is the usual Laplacian and $\mathbf{1}_{\mathbb S}$ is the indicator function of the set $\mathbb S$.

In this article, we take $\beta=2$. As a consequence, the point process given in \eqref{eq:2D CG} is determinantal with a kernel $K_n$ (see \eqref{eq: def kernel} for a precise definition). Further, we consider a model where the potential is radially symmetric and is such that a spectral gap is about to appear in the bulk around a circle of radius $r_*>0$ and thus the number of components of the droplet increases. A more precise description of the model is given in \autoref{pres of model}. In this setting, we derive a large $n$ expansion of the partition function $Z_n^{(2)}$. This is of particular interest, as the asymptotic expansion of the partition function captures important information about the system, such as the logarithmic energy, the entropy, some geometric properties and more. 
 In \cite{ZW06}, Zabrodin and Wiegmann predicted that, for sufficiently regular potentials---namely, those with smooth boundaries and no singularities or critical behavior---the partition function $Z_n^{(\beta)}$ admits an asymptotic expansion of the following form, as $n\rightarrow \infty$,
\begin{equation}
    \log( Z_n^{(\beta)})=\mathcal{C}_0 n^2+\mathcal{C}_1 n\log(n)+\mathcal{C}_2n+\mathcal{C}_3\log(n) +\mathcal{C}_4+O(n^{-1}).
\end{equation}
and proposed explicit formulas for the constant $\mathcal{C}_j \equiv \mathcal{C}_j(\beta,Q)$, $(j=0,\dots, 4)$. Recently, it was shown in \cite{HM12,DGZ14} and \cite{LS18} that, as $n\rightarrow\infty$,

\begin{equation}
    \log( Z_n^{(\beta)})=-\frac{\beta}{2}n^2 I_Q[\mu_Q]+\frac{\beta}{4} n\log(n)-n\left[\mathcal{C}(\beta)+\left(1-\frac{\beta}{4}\right)E_Q[\mu_Q]\right] +o(n).
\end{equation}
where $\mathcal{C}(\beta)$ is independent of $Q$ and
\begin{equation}\label{eq:entropy}
    E_{Q}[\mu_Q]:=\int_{\mathbb C}\log (\frac{\Delta Q}{4})d\mu_Q,
\end{equation}
is the entropy associated with $\mu_Q$. The quantities $\mu_Q$ and $I_Q[\mu_Q]$ are given in \eqref{eq:eqmeas} and \eqref{eq: energy}, respectively. In the determinantal setting this conjecture has been partially verified in the works of \cite{ByunKangSeoYang2025FreeEnergySpherical, B25,BKS23,ACC23}.

On the real line, the vanishing of the equilibrium measure in the interior and multiple cuts have been investigated. In particular, we refer to \cite{B09,CL07,ClaeysGravaMcLaughlin2015Asymptotics,CFWW25} for further details in this direction.

In two dimensions, the statistical properties near the bulk and the outer boundary of the droplet are well studied; see, for example, \cite{AC22,HW21,FH99,FJ96,LR16,AKM18,S22,CW26}. While the study of spectral gaps is an emerging area of focus in the context of Coulomb gases, the presence of a spectral gap with hard edges has been analysed in \cite{C24,ACCL24,B25}, where one observes the emergence of oscillatory terms involving theta functions in the constant term of the free energy expansion. See also \cite{BG24,CFWW25} for an analogous result in one dimension.

For soft edges, the authors in \cite{AKS21} consider the emergence of a spectral gap at the origin and demonstrate that the kernel can be expressed in terms of Mittag-Leffler functions $E_{a,b}$. Furthermore, in \cite{ACC23,ACC22a}, the free energy and fluctuations of linear statistics are studied in the critical regime, where a new component of the droplet emerges away from the primary droplet, a phenomenon referred to as the existence of a
``shallow outpost''. They show that the number of particles near the outpost follows a specific Heine distribution. Recently, the question of universality in this setting has been addressed in \cite{AC24}. See also the recent work of \cite{K26} for multiple outposts. 

\begin{figure}[h!]
        \centering
        \includegraphics[width=0.45\textwidth]{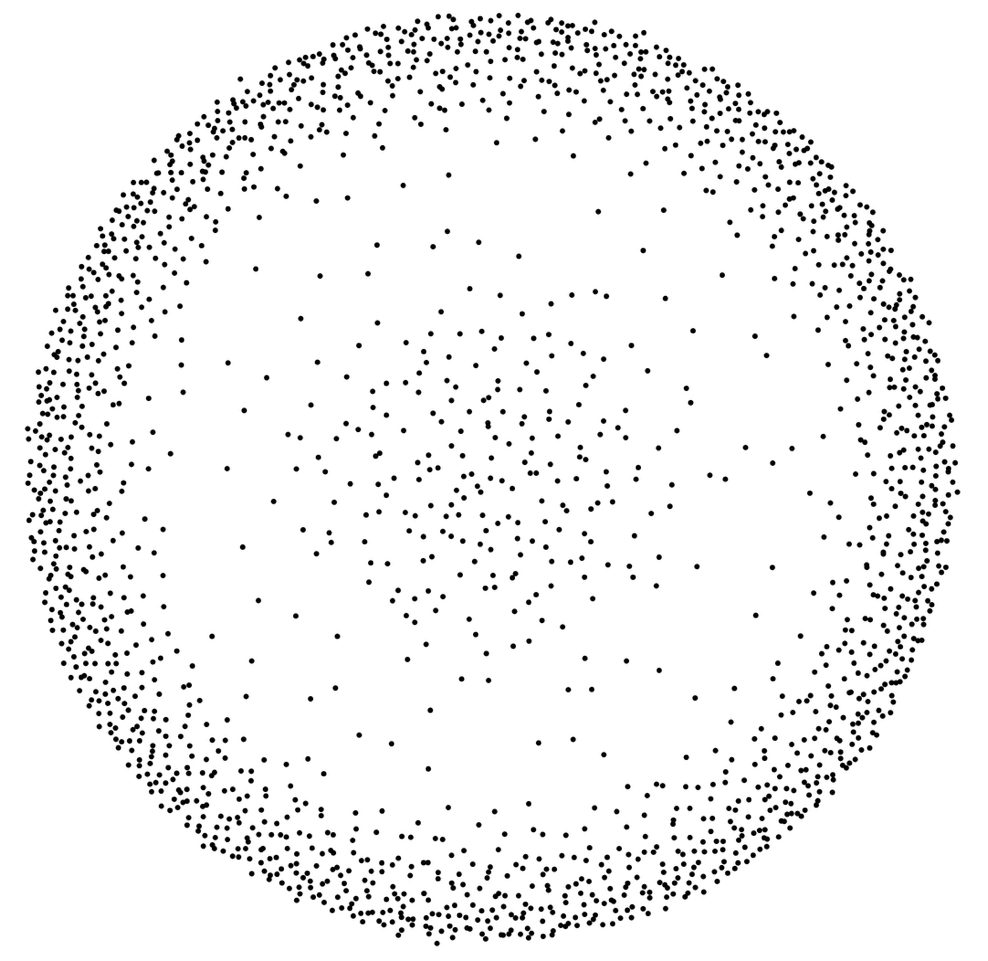}
        \vspace{5pt} 
        \caption{2D Coulomb gas droplet with a nearly formed circular gap. Illustration generated using the Metropolis-adjusted Langevin algorithm for the potential $Q$ given in \eqref{eq: example potential}, courtesy of Aron Wennman. }  
        \label{fig: example CG}
\end{figure}
In contrast to \cite{AKS21}, we investigate, in the present work, the emergence of a spectral gap away from the origin, localized around a one-dimensional circle, as illustrated in \autoref{fig: example CG}. Our result for the expansion of the partition function features a novel term of order $n^{1/4}$ (see \autoref{theo1}), which is new in the literature, as well as a more general integral formulation of the constant terms (see \eqref{eq: functional reduce}) compared with \cite[Theorems 1.1 and 1.2]{BKS23}. This scaling $n^{1/4}$ is inversely proportional to the mean level spacing of the particles in the vicinity of the circle of radius $r_*$. 

If one zooms precisely at this scale around $r_*$ in the radial direction, the local correlations between particles are governed by a new correlation kernel $K_*$ \eqref{eq: K_*} (see \autoref{rem: local PP})  which differs from the Ginibre kernel governing local correlations in the bulk; we refer to \cite{AKM18,AC22} for a precise definition. The new kernel is expressed in terms of a Pearcey-like integral \cite{JP16}; see \autoref{theo2}. 

Before presenting the main results, let us give a more precise description of the model under consideration.

\section{Main results}
\subsection{Presentation of the model}\label{pres of model}
\begin{figure}[h!]
        \centering
        \includegraphics[width=0.8\textwidth]{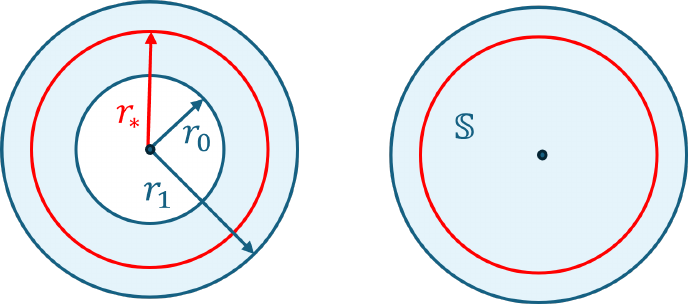}
        \vspace{5pt} 
        \caption{\textit{(Left)} $r_0>0$: The droplet $\mathbb{S}$ (blue) is an annulus $\mathbb{A}_{r_0,r_1}$ and $\mathbb{T}_{r_*}$, the circle of radius $r_*$ (red), lies inside the annulus. \textit{(Right)}  $r_0=0$: The droplet $\mathbb{S}$ (blue) is a disk of radius $r_1$, $\mathbb{D}_{r_1}:=\mathbb{A}_{0,r_1}$, and $\mathbb{T}_{r_*}$ (red), lies inside the disk.}  
        \label{fig: config}
\end{figure}

We consider here the distribution \eqref{eq:2D CG} with a smooth and radially symmetric potential $Q\in\mathrm{C}^6(\mathbb{S})$. We denote
\Eq{eq:complex potential}{
Q(z):=q(\abs{z}).
}
To impose the critical condition corresponding to the near emergence of a gap, we require that the Laplacian of $Q$ vanishes on the circle $\mathbb{T}_{r_*}:=\{z\in \mathbb{C}\mid \abs{z}=r_*\}$ and to be positive everywhere else, i.e.
\Eq{eq: subharm cond}{
\Delta Q(r_*)=0,\quad \Delta Q(r)>0,\quad r\neq r_*>0.
}
As $Q$ is smooth, this ensures that $r_*$ is the global minimum of $r\mapsto \Delta Q(r)$. The radius $r_*$ is thus a double zero of the function $r\mapsto \Delta Q(r)$ and we impose that it is, however, not a triple zero, i.e. 
\begin{equation}\label{eq: vanish laplacian}
    \partial_r \Delta Q(r_*)=0, \qquad  \partial^2_r\Delta Q(r_*)>0.
\end{equation}
This latter condition is crucial as it determines the criticality of the formation of gap  phenomenon; see \autoref{Statement of results} $\S$ below \autoref{cor:mls}. Additionally, we assume the potential $Q$ to have sufficient growth at infinity,
\Eq{eq: growth pot}{
\liminf_{\abs{z}\to \infty}\frac{Q(z)}{2\log(\abs{z})}>1,
}
which guarantees $\mu_Q$ to be compactly supported. The weak subharmonicity condition \eqref{eq: subharm cond} ensures that the support of the equilibrium measure $\mu_Q$ is
\Eq{}{
\mathbb{S}:=\mathrm{supp}(\mu_Q)=\mathbb{A}_{r_0,r_1},\quad 0\leq r_0<r_*<r_1,
}
where $\mathbb{A}_{r_0,r_1}:=\{z\in \mathbb{C}\mid r_0\leq \abs{z}\leq r_1 \}$. \autoref{fig: config} illustrates the configuration and a simple example of potential verifying all the previous conditions is presented in \autoref{example potential}. 

We will further restrict ourselves to the case $\beta=2$. Thus, the density of the point process can be rewritten in the determinantal form
\Eq{eq: dpp}{
\frac{1}{Z_n^{(2)}}\prod_{1\leq j<k\leq n}|z_j-z_k|^{2}\prod_{j=1}^{n}e^{-n\,Q(z_j)}=\frac{1}{n!}\det\left[K_n(z_j,z_k) \right]_{j,k=1}^n,
}
where the kernel $K_n$ can be explicitly expressed in terms of the weighted monic orthogonal polynomials
\Eq{}{
p_{j,n}(z):=z^j e^{-\frac{n}{2}Q(z)},\quad j\in\{ 0,1,2,\dots, n-1\},
}
in the following way
\begin{equation}\label{eq: def kernel}
K_{n}(z_1,z_2)=\sum_{j=0}^{n-1}\frac{p_{j,n}(z_1)\overline{p_{j,n}(z_2)}}{{\norm{p_{j,n}}^2}},
\end{equation}
with the norm of the polynomials given by
\begin{equation}\label{eq: norms}
\norm{p_{j,n}}^2=\int_{\mathbb{C}} \abs{z}^{2j} e^{-nQ(\abs{z})} \frac{dA (z)}{\pi}=2\int_0^\infty  r^{2j+1} e^{-nq(r)}dr.
\end{equation}
\rem{}{
We recall that, in general, the kernel of a determinantal point process is not uniquely defined. Indeed, for any non-vanishing function $g$, known as a cocycle, the kernel given by $\left(g(z_1)/g(z_2) \right) K_n(z_1, z_2)$ defines the same point process. This invariance can be verified through elementary row and column operations on the determinant.
}
The partition function $Z_n^{(2)}$ then reads 
\begin{equation}\label{eq:betapart}
   Z_n:=Z_n^{(2)}=\int_{\mathbb C^n}\prod_{1\leq j<k\leq n}|z_j-z_k|^{2}\prod_{j=1}^{n}e^{-n Q(z_j)}dA(z_j).
\end{equation}
A standard calculation, as shown in \cite[Exercises 15.3, q.1(ii)]{F10}, shows that the corresponding $n$-dimensional integral $Z_n$ factorizes into a product of one-dimensional integrals, which, after taking the logarithm, yields
\Eq{eq:defuj}{
\log(\frac{Z_n}{n!\,(2\pi)^n })=\sum_{j=0}^{n-1}\log(u_{j}), \qquad u_j:=\int_{0}^{\infty}r^{2j+1}e^{-n q(r)}dr.
}

Before stating our main results, we introduce the following constants
\Eq{eq:vcon}{
\gamma_Q &:= \partial^2_r \Delta Q(r_*)= q^{(4)}(r_*) + \frac{6 q'(r_*)}{r_*^3}, &\qquad \varkappa_Q &:=\frac{r_*}{2}\left(\partial_r^2\Delta Q(r_*)\right)^{1/4}
}
as well as the following Pearcey-like integral 
\Eq{eq: pearcey}{
\mathcal{P}(\theta):=\int_{-\infty}^{\infty} e^{-\theta P(y)}dy,\qquad P(y):=\frac{y^2}{4}+ \frac{y^3}{3!}+\frac{y^4}{4!},
}
which appears both in the large $n$ expansion of the partition function (cf. \autoref{theo1}) and in the expression of the local correlation kernel (cf. \autoref{theo2} and \autoref{coro2}).

\subsection{Statement of results}\label{Statement of results}
The following theorem presents our first main result: a large-$n$ expansion of the partition function $Z_n$, which features a novel term of order $n^{1/4}$, as well as a constant term that differs from the non-critical case studied in \cite[Theorems 1.1 and 1.2]{BKS23}. For each configuration---annulus and disk (see \autoref{fig: config})---the constant term is expressed in terms of the following functionals:
\Eq{eq: def FQ}{
\mathcal{F}_{Q}[\mathbb{A}_{r_0,r_1}]:=&\int_{r_0}^{r_1}\left[-\frac{7}{24}r \left(\frac{\partial_r \Delta Q(r)}{\Delta Q(r)}\right)^2
+\frac{1}{3 } r\frac{\partial_r^2 \Delta Q(r)}{\Delta Q(r)}
-\frac{r C_{2,Q}+C_{1,Q}}{(r-r_*)^2}\right]dr\\
&- \frac{19}{48 }\left[ r_1\frac{\partial_r \Delta Q(r_1)}{\Delta Q(r_1)}-r_0\frac{\partial_r \Delta Q(r_0)}{\Delta Q(r_0)}\right]-\frac{1}{6}\log(\frac{r_1}{r_0})+\frac{1}{4}\log\!\left(\frac{\Delta Q(r_1)}{\Delta Q(r_0)}\right)\\
&-\left[\frac{r_1 C_{2,Q}+C_{1,Q}}{r_1-r_*}
+\frac{r_0 C_{2,Q}+C_{1,Q}}{r_*-r_0}\right]+C_{2,Q} \log\!\left(\abs{\frac{r_1-r_*}{r_*-r_0}}\right),
}
with 
\Eq{eq: def CQ}{
C_{1,Q}:=\frac{7r_*}{24},\qquad C_{2,Q}:=\frac{1}{2}\left(
\frac{1}{9}\frac{r_* \partial_r^3\Delta Q(r_*)}{ \partial_r^2\Delta Q(r_*)} - 1\right),
} 
and 
\Eq{eq: def FQ2}{
\mathcal{F}_{Q}[\mathbb{D}_{r_1}]:=&\int_{0}^{r_1}\left[-\frac{7}{24}r \left(\frac{\partial_r \Delta Q(r)}{\Delta Q(r)}\right)^2
+\frac{1}{3 } r\frac{\partial_r^2 \Delta Q(r)}{\Delta Q(r)}
-\frac{r C_{2,Q}+C_{1,Q}}{(r-r_*)^2}\right]dr\\
&- \frac{19}{48 }r_1\frac{\partial_r \Delta Q(r_1)}{\Delta Q(r_1)}-\frac{1}{6}\log(r_1)+\frac{1}{4}\log\!\left(\frac{\Delta Q(r_1)}{\Delta Q(0)}\right)-\left[\frac{r_1 C_{2,Q}+C_{1,Q}}{r_1-r_*}+\frac{C_{1,Q}}{r_*}\right]\\
&+C_{2,Q} \log\!\left(\abs{\frac{r_1-r_*}{r_*}}\right)-\frac{1}{12 }\log\!\left(\frac{\Delta Q(0)}{4}\right).
}
We note that both functionals, when applied to a non-critical potential $\hat{Q}$---for which $\Delta \hat{Q}>0$ on the full droplet $\mathbb{S}$---reduce, respectively, to the functionals $F_{\hat{Q}}[\mathbb{A}_{r_0,r_1}]$ and $F_{\hat{Q}}[\mathbb{D}_{r_1}]$ defined in \cite[Theorems 1.1 and 1.2]{BKS23} (where $\Delta$ is defined as one quarter of the usual Laplacian). Indeed, for non-critical potentials, $r_*$ either lies outside $[r_0, r_1]$ or is undefined---in the latter case, one may take $r_* > r_1$; then, integration by parts on \eqref{eq: def FQ} and \eqref{eq: def FQ2} yields
\Eq{eq: functional reduce}{
\mathcal{F}_{\hat{Q}}[\mathbb{A}_{r_0,r_1}]=F_{\hat{Q}}[\mathbb{A}_{r_0,r_1}],\qquad\mathcal{F}_{\hat{Q}}[\mathbb{D}_{r_1}]=F_{\hat{Q}}[\mathbb{D}_{r_1}].
}
This is expected, as in this case one is working within the framework considered by \cite{BKS23}.
\theo{theo1}{
Let the potential $Q$ be such as in \autoref{pres of model}. \begin{itemize}
    \item \textbf{Annulus case} ($r_0>0$), see \autoref{fig: config}. The large $n$ expansion of the partition function $Z_n$ \eqref{eq:betapart} is given as 
 \Eq{}{
      \log\left( \frac{Z_n}{n!\,(2\pi)^n}\right)=&\,-n^2 I_Q[\mu_Q]-\frac{1}{2}n\log(n)+ n\left[\frac{ \log(\pi/2) }{2}-\frac{1}{2}E_Q[\mu_Q]\right]\\
      &+n^{1/4}\frac{\varkappa_Q}{2}\int_{-\infty}^{\infty}x^2\log(\frac{x^2}{\sqrt{4\pi}} \mathcal{P}(x^4))dx+\mathcal{F}_{Q}[\mathbb{A}_{r_0,r_1}]+o(1).
 }
 \item \textbf{Disk case} ($r_0=0$), see \autoref{fig: config}. The large $n$ expansion of the partition function $Z_n$ \eqref{eq:betapart} is given as 
 \Eq{}{
      \log\left( \frac{Z_n}{n!\,(2\pi)^n}\right)=&\,-n^2 I_Q[\mu_Q]-\frac{1}{2}n\log(n)+ n\left[\frac{ \log(\pi/2) }{2}-\frac{1}{2}E_Q[\mu_Q]\right]\\
      &+n^{1/4}\frac{\varkappa_Q}{2}\int_{-\infty}^{\infty}x^2\log(\frac{x^2}{\sqrt{4\pi}} \mathcal{P}(x^4))dx\\
      &-\frac{1}{12}\log(n)+\zeta'(-1)+\mathcal{F}_{Q}[\mathbb{D}_{r_1}]+o(1).
 }
 \end{itemize}
The definitions of $I_Q[\mu_Q]$, $E_Q[\mu_Q]$, $\mathcal{F}_{Q}[\mathbb{A}_{r_0,r_1}]$, $\mathcal{F}_{Q}[\mathbb{D}_{r_1}]$, $\varkappa_Q$ and $\mathcal{P}$ are given, respectively, in \eqref{eq: energy}, \eqref{eq:entropy}, \eqref{eq: def FQ}, \eqref{eq: def FQ2}, \eqref{eq:vcon} and \eqref{eq: pearcey}. Here, $\zeta$ is the Riemann zeta function.
}
\rem{}{
We emphasize that the only structural difference with \cite[Theorems 1.1 and 1.2]{BKS23} is the presence of the term of order $n^{1/4}$. As discussed above, the integral formulation of the constant terms is merely a generalization of those in \cite{BKS23}, accounting for the singularities of the integrand. The term of order $\log(n)$ remains unchanged, as it depends on the Euler characteristic of the droplet or the presence of a hard wall (cf. \cite{AFLS25}). In our setting, there is no hard wall, and we consider droplets with the same Euler characteristic as in \cite{BKS23}, namely, either a disk or an annulus.
}
It is interesting to notice that the novel term of order $n^{1/4}$ comes with a universal constant---independent of the potential---given by the convergent integral (see \eqref{eq:asy})
\Eq{eq:unicon}{
\int_{-\infty}^{\infty}x^2\log(\frac{x^2}{\sqrt{4\pi}} \mathcal{P}(x^4))dx=\int_0^\infty \sqrt{y}\log(\frac{y}{\sqrt{4\pi}} \mathcal{P}(y^2))dy\approx 1.603,
}
which involves the Pearcey-like integral $\mathcal{P}$ \eqref{eq: pearcey}. The physical significance of this term is likely related to surface tension. 
In the case of a 2D Coulomb gas in a confining potential (see \cite{Jancovici1994,AFLS25}), there exists a surface tension contribution of order $n^{1/2}$, whose sign aligns with that of the leading terms of order $n^2$ and $n$. This contrasts with the present setting, where the contribution appears at a lower fractional order and carries the opposite sign. The fractional order reflects the degree of criticality of the potential and will be addressed in detail below.

\begin{figure}[h!]
        \centering
        \includegraphics[width=\textwidth]{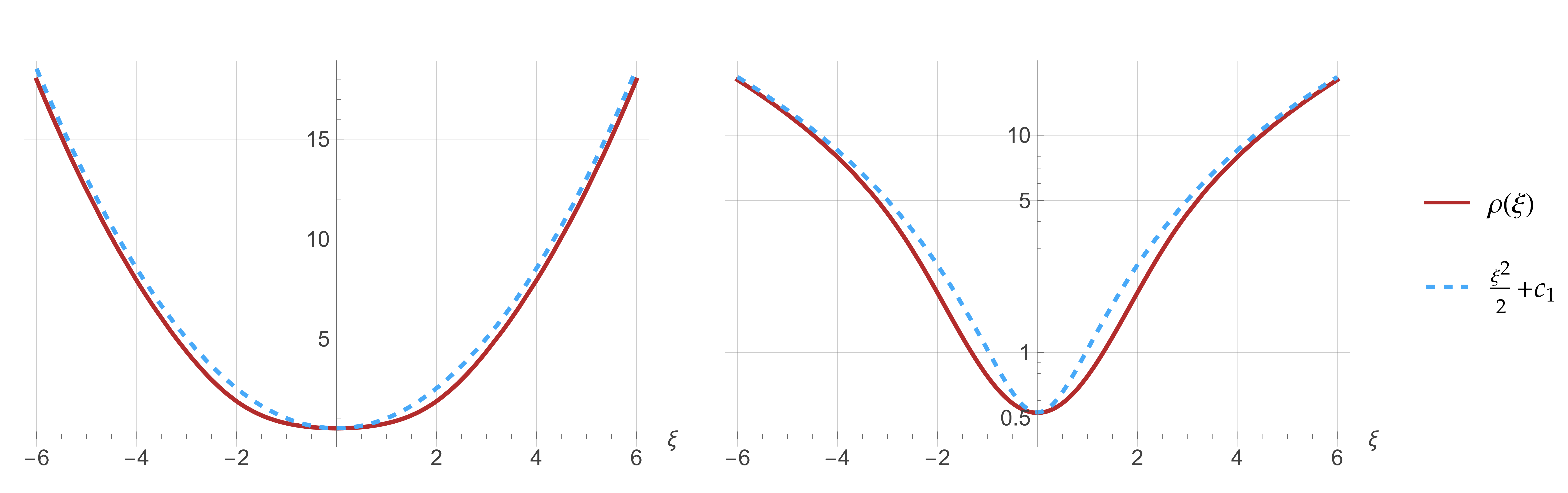}
        \caption{Plots of the functions $\xi \mapsto \rho(\xi)$, as defined in equation~\eqref{eq: def rho}, and $\xi \mapsto \xi^2/2 + c_1$, where $c_1 = \rho(0)$, shown for comparison over the real interval $\xi \in [-6, 6]$. \textit{(left)} normal scale. \textit{(right)} $y$-axis in $\log$ scale. }  
        \label{fig: 1pt fct1}
\end{figure}

\begin{figure}[h!]
        \centering
        \includegraphics[width=0.8\textwidth]{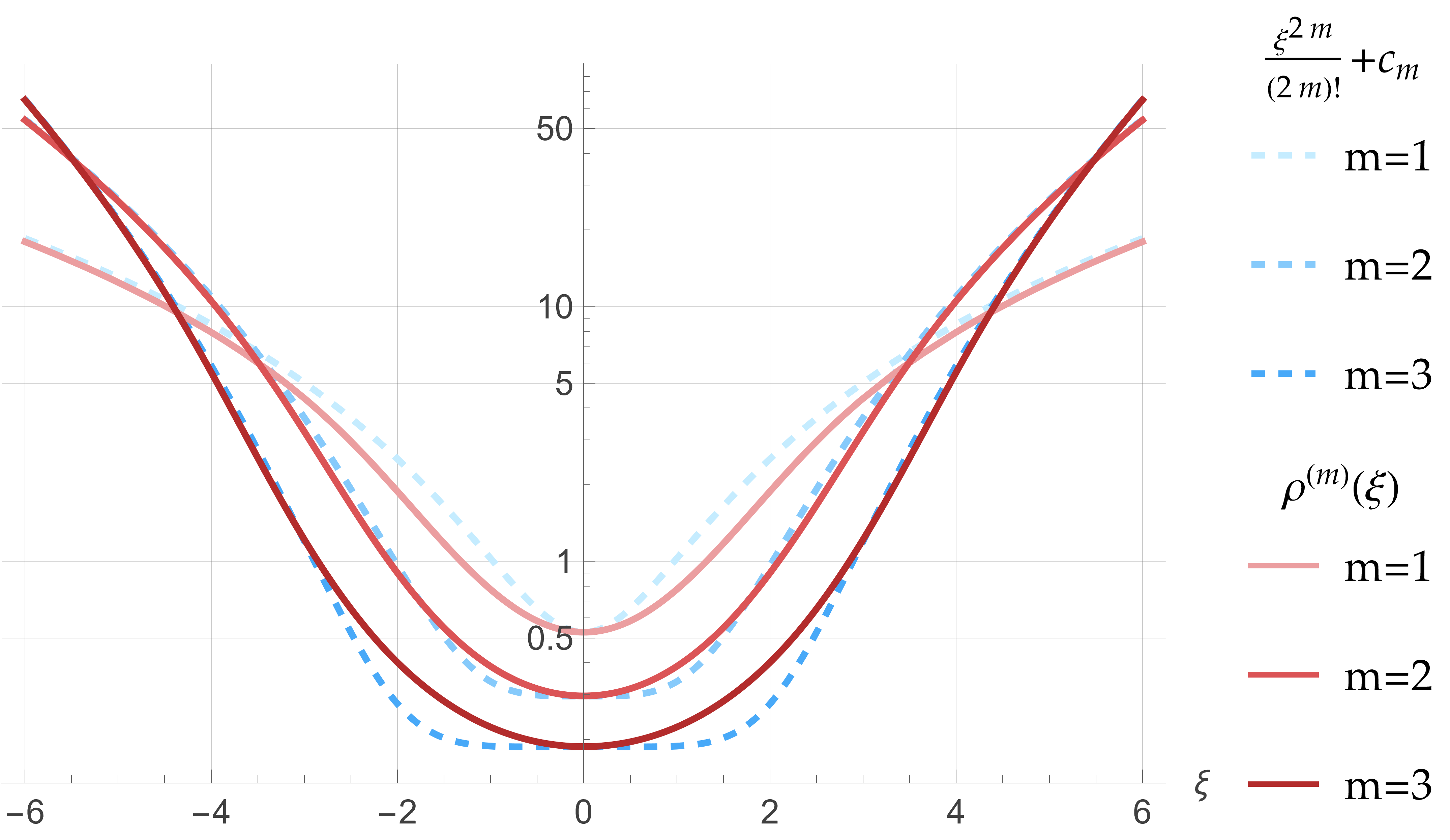}
        \caption{Plots of the functions $\xi \mapsto\rho^{(m)}(\xi)$, defined in \eqref{eq: def rho m}, and $\xi\mapsto \xi^{2m}/(2m)!+c_m$, with $c_m=\rho^{(m)}(0)$, shown for comparison over the real interval $\xi \in [-6, 6]$. For $m=1,2,3$. $y$-axis in $\log$ scale. }  
        \label{fig: rho m}
\end{figure}

Turning to the local correlations between particles in the vicinity of the circle of radius $r_*$, the second main result is the double scaling limit of the kernel and is presented in the following theorem. We have obtained an explicit expression for the first correction term of order $n^{-1/4}$. However, since it is not universal (cf. \autoref{rem: univ correct}), cumbersome, and offers little additional insight, we present it in \autoref{app correction term} for completeness.
\theo{theo2}{
   For $k=1,2$, let $z_k=r_*+\frac{\xi_k}{(\gamma_Q\, n)^{1/4}}\in \mathbb C$, with $\gamma_Q$ given in \eqref{eq:vcon} and $|\xi_k|\leq \log(n)$. Under the same assumptions as \autoref{theo1}, the kernel \eqref{eq: def kernel}, up to a cocycle, admits the following uniform expansion, as $n \to \infty$,
    \begin{equation}
    \begin{aligned}\label{eq: double scaling limit}
       \frac{1}{\sqrt{\gamma_Q\,n}} K_{n}(z_1,z_2)
       =&\,\frac{1}{8}\int_{-\infty}^{\infty}e^{i\Im(\xi_1-\xi_2) \frac{x^3}{12}}\frac{\exp\left(-\frac{x^4}{2}\left[ P\left(\frac{\Re(\xi_1)}{x}-1\right)+ P\left(\frac{\Re(\xi_2)}{x}-1\right)\right]\right)}{\mathcal{P}(x^4)}\abs{x}dx\\&+O(n^{-1/4}),
  \end{aligned}
    \end{equation}
where $\mathcal{P}$ and $P$ are given in \eqref{eq: pearcey}.
}
\rem{}{Due to radial symmetry of the problem, the result of \autoref{theo2} holds for any $\tilde{z}_k=e^{i\varphi}z_k$, $\varphi\in\mathbb R$. We choose $\varphi=0$ for simplicity. Let us also point out that another choice of local coordinates is possible, namely,
\Eq{}{
z_k=\left(r_*+\frac{\nu_k}{(\gamma_Q\, n)^{1/4}}\right)e^{i \frac{\vartheta_k}{r_*(\gamma_Q\,n)^{1/4}}},\quad \nu_k,\vartheta_k\in\mathbb{R},
}
as it is done in \cite{ACC22a}. The two choices are related via the expansion
\Eq{}{
z_k=r_*+\frac{\xi_k}{(\gamma_Q\, n)^{1/4}}=\left(r_*+\frac{\Re(\xi_k)}{(\gamma_Q\, n)^{1/4}}\right)e^{i \frac{\Im(\xi_k)}{r_*(\gamma_Q\,n)^{1/4}}}+O(n^{-1/2}),
}
where one can set $\nu_k=\Re(\xi_k)$ and $\vartheta_k=\Im(\xi_k)$.
}
The above theorem makes clear that the sequence of kernels $K_n$ admits a double scaling limit---being the leading term of the expansion \eqref{eq: double scaling limit}---when considering a pair of points whose radial distance from $r_*$ is of order $n^{-1/4}$. This is to be expected, as this scale corresponds to the mean level spacing of particles in the vicinity of $r_*$; cf. \autoref{cor:mls}. Moreover, this limiting kernel is proportional to a universal kernel---independent of the potential---which we denote
\Eq{eq: K_*}{
K_*(\xi_1,\xi_2):=\frac{1}{2}\int_{-\infty}^{\infty}e^{i\Im(\xi_1-\xi_2) \frac{x^3}{12}}\frac{\exp\left(-\frac{x^4}{2}\left[ P\left(\frac{\Re(\xi_1)}{x}-1\right)+ P\left(\frac{\Re(\xi_2)}{x}-1\right)\right]\right)}{\mathcal{P}(x^4)}\abs{x}dx.
}
\begin{figure}[h!]
        \centering
        \includegraphics[width=0.6\textwidth]{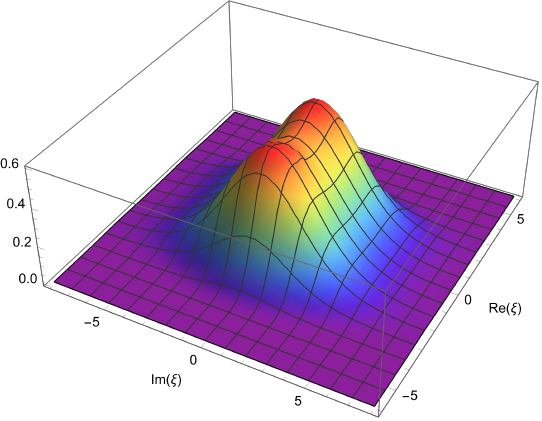}
        \caption{Plot of the function $\xi\mapsto\abs{K_*(\xi,0)}$, as defined in equation~\eqref{eq: K_*} for $-6\leq \Re(\xi)\leq6$,  $-8\leq \Im(\xi)\leq8$. The quantity $\abs{K_*(\xi,0)}^2$ is a measure of the covariance between a point at $z=0$ and a point at $z=\xi$ (in the local referential centred at $r_*$). See \autoref{rem: interpretation} for some further comments.}  
        \label{fig: kernel3D}
\end{figure}
A plot of the function $\xi\mapsto\abs{K_*(\xi,0)}$ is given in \autoref{fig: kernel3D} where one can notice the symmetry $K_*(\overline{\xi},0)=K_*(-\xi,0)$. Along with the more trivial symmetries $\abs{K_*(\xi,0)}=\abs{K_*(0,\xi)}=\abs{K_*(0,\overline{\xi})}=\abs{\overline{K_*(0,\xi)}}$, which are only true when taking the modulus.
\rem{rem: local PP}{
The limiting kernel $K_*$ \eqref{eq: K_*} defines a new local point process for which the $k$-point correlation functions $R_k$ are then given by 
\Eq{}{
R_k(\xi_1,\ldots,\xi_k)=\det\left[K_*(\xi_i,\xi_j)\right]_{i,j=1}^k.
}
It is important to note that, unlike the initial macroscopic point process, the microscopic $k$-point correlation functions cannot be normalized to probability densities. This justifies our omission of the overall factor $1/4$ in the definition of $K_*$. However, we retain a factor of $1/2$---reflecting the critical nature of the potential---to allow for meaningful comparison between the limiting one-point functions at different criticalities; see \eqref{eq: fct m} and \autoref{fig: rho m}.
}
Setting $z_1 = z_2$ in \autoref{theo2} immediately yields the local $1$-point correlation function as a corollary. The explicit expression for the correction term of order $n^{-1/4}$ is obtained from \autoref{theo2 strong}.
\coro{coro2}{
Let $z=r_*+\frac{\xi}{(\gamma_Q\, n)^{1/4}}\in\mathbb C$, with $|\xi|\leq \log(n)$. Then, as $n\rightarrow\infty$, we have 
\begin{equation}
       \frac{1}{\sqrt{\gamma_Q\,n}} K_{n}(z,z)=\frac{1}{8}\int_{-\infty}^{\infty}\frac{\exp\left(-x^4 P\left(\frac{\Re(\xi)}{x}-1\right)\right)}{\mathcal{P}(x^4)}\abs{x}dx+O(n^{-1/4})
    \end{equation}
where $\mathcal{P}$ and $P$ are given in \eqref{eq: pearcey} and  $\gamma_Q$ given in \eqref{eq:vcon}.
}
The local $1$-point correlation function is then proportional to the universal positive function
\Eq{eq: def rho}{
\rho(\xi):=K_*(\xi,\xi)=\frac{1}{2}\int_{-\infty}^{\infty}\frac{\exp\left(-x^4 P\left(\frac{\Re(\xi)}{x}-1\right)\right)}{\mathcal{P}(x^4)}\abs{x}dx,
}
for which a straightforward application of Laplace's method using the change of variable $x\mapsto \Re(\xi)/(1+x)$ and the asymptotic of $\mathcal{P}$ \eqref{eq:asy} yields the asymptotic equivalence 
\Eq{eq: asymp rho}{
\rho(\xi)\sim \frac{\Re(\xi)^{2}}{2},\quad \Re(\xi)\to\infty.
}
Since $\rho$ depends only on the real part of its argument, in \autoref{fig: 1pt fct1} we plot $\rho$ alongside the comparison function $\xi \mapsto \xi^2/2 + \rho(0)$ over a real interval centered at $0$.
\rem{rem: interpretation}{
Comparing \autoref{fig: kernel3D} with \autoref{fig: 1pt fct1} one can notice that the local 2-point correlation function, which is proportional to 
\Eq{}{
\det 
\begin{bmatrix}
K_*(\xi, \xi) & K_*(0, \xi) \\
K_*(\xi, 0) & K_*(0, 0)
\end{bmatrix}=\rho(\xi)\rho(0)-\abs{K_*(\xi, 0)}^2,
}
is dominated by the product $\rho(\xi)\rho(0)$ for $\abs{\xi}$ a bit larger than $2$. Thus, not only does $\abs{K_*(\xi, 0)}^2$ drop very fast for $\abs{\xi}$ larger than $2$ but it is also completely negligible compared to the product $\rho(\xi)\rho(0)$. As $\abs{K_*(\xi, 0)}^2$ is a measurement of the covariance between a point at $z=0$ and $z=\xi$, this means that beyond a distance of $\abs{\xi}\approx 2$, the points are more or less statistically independent, i.e. the particles do not feel each other. The two maxima that one can see in \autoref{fig: kernel3D} and located around $\Re(\xi)\approx \pm 1$ correspond to the average location where two particles are the most correlated, i.e. where the particles will feel each other the most.
}

It is straightforward to verify that $\rho$ is an even function, as a simple change of variable $x\mapsto-x$ in the integral yields $\rho(-\xi)$. This is heuristically expected since, in the vicinity of $r_*$, particles have no way to distinguish between the bulk above $r_*$ and the bulk below $r_*$, therefore, the probability for a particle to be on either side should be the same.

An immediate consequence of \autoref{theo2}---visible from the appropriate scaling of $z$---is the scale of the mean level spacing of particles near $r_*$. This quantity is defined as the average distance from the location $z = r_*$ to the nearest point of the process. A characterization, which we will in fact adopt as the precise mathematical definition, is given below.
\defi{def: mls}{
The \emph{local mean level spacing} at $r_*$ is the radius $s_n$ determined by 
\begin{equation}\label{def:mls}
    \mathbb{E}\left[N_n\!\bigl(\mathbb{D}(r_*,s_n)\bigr)\right]=1,
\end{equation}
where $\mathbb{D}(r_*, d):=\{z\in\mathbb{C}\mid \abs{z-r_*}<d\}$ and $N_n(\Omega)$ denotes the number of points in a given realization of the point process that fall within the set $\Omega$.
}

\coro{cor:mls}{For the determinantal point process with kernel $K_n$ given in \autoref{theo2} with potential $Q$ such as in \autoref{pres of model}, the mean level spacing of the particles around a point $z\in\mathbb{C}$ with $|z|=r_*$ is proportional to $n^{-1/4}$.}
The particular scaling $n^{-1/4}$ is a direct consequence of the degree of criticality of the potential; see \autoref{pres of model}. More precisely, if one considers a potential such that 
\Eq{eq: m vanish deriv}{
\partial_r^{j} \Delta Q(r_*)=0, \quad j=0,\ldots,2m-1, \qquad  \partial^{2m}_r\Delta Q(r_*)>0,
}
for an integer $m\geq 1$, then, the mean level spacing is expected to be proportional to $n^{-\frac{1}{2m+2}}$. Note that the first nonvanishing derivative has to be even for $r_*$ to be an extremum, otherwise it is simply an inflection point. For $m>1$, one expects the universal constants to be
\Eq{eq: cst m}{
\gamma_{Q,m} &:= \partial_r^{2m} \Delta Q(r_*), &\qquad \varkappa_{Q,m} &:=\frac{r_*}{2 }\left(\partial_r^{2m}\Delta Q(r_*)\right)^{\frac{1}{2m+2}},
}
and the auxiliary functions to be
\Eq{eq: fct m}{
\mathcal{P}_{m}(\theta):=\int_{-\infty}^{\infty} e^{-\theta P_{m}(y)}dy,\qquad P_m(y):=\sum_{k=2}^{2m+2} \frac{y^k}{(2m+2-k)!\,k!}.
}
This leads to the limiting kernel
\Eq{eq: kernel m}{
K_*^{(m)}(\xi_1,\xi_2):=\frac{1}{(2m)!}&\int_{-\infty}^{\infty}e^{i\Im(\xi_1-\xi_2) \frac{x^{2m+1}}{2(2m+1)!}}\\
&\times \frac{\exp\left(-\frac{x^{2m+2}}{2}\left[P_{m}\left(\frac{\Re(\xi_1)}{x}-1\right)+P_{m}\left(\frac{\Re(\xi_2)}{x}-1\right)\right]\right)}{\mathcal{P}_{m}(x^{2m+2})}\abs{x}^{2m-1}dx.
}
\rem{}{
While we express the limiting kernel in terms of the function $P_m$ to make the structure more apparent, we note that one can use the binomial identity to rewrite it as a homogeneous polynomial 
\Eq{}{
x^{2m+2}P_{m}\left(\frac{y}{x}-1\right)=\frac{1}{(2m+2)!}\left(y^{2m+2}-(2m+2)y\,  x^{2m+1}+(2m+1)x^{2m+2} \right).
}
}
Gathering the above equations gives rise to the following conjectures.
\conj{conj: Z_n}{
Let the potential $Q$ be as in \autoref{pres of model}, now assuming $Q \in \mathrm{C}^{2m+4}(\mathbb{S})$ and satisfying \eqref{eq: m vanish deriv} in place of \eqref{eq: vanish laplacian}, for some integer $m \geq 2$. In both the annulus case ($r_0>0$) and the disk case ($r_0=0$) (see \autoref{fig: config}), the large $n$ expansion of the partition function $Z_n$ \eqref{eq:betapart} is given as 
 \Eq{}{
      \log\left( \frac{Z_n}{n!\,(2\pi)^n}\right)=&\,-n^2 I_Q[\mu_Q]-\frac{1}{2}n\log(n)+ n\left[\frac{ \log(\pi/2) }{2}-\frac{1}{2}E_Q[\mu_Q]\right]\\
      &+n^{\frac{1}{2m+2}}\frac{\varkappa_{Q,m}}{(2m)!}\int_{-\infty}^{\infty}x^{2m}\log(\frac{\abs{x}^{m+1}}{\sqrt{2\pi (2m)!}} \mathcal{P}_m\left(x^{2m+2}\right))dx+O(\log(n)),
 }
where $I_Q[\mu_Q]$, $E_Q[\mu_Q]$, $\varkappa_{Q,m}$ and $\mathcal{P}_m$ are given, respectively, in \eqref{eq: energy}, \eqref{eq:entropy}, \eqref{eq: cst m} and \eqref{eq: fct m}.
}
\rem{}{
Regarding the sub-leading terms, they are expected to depend on whether the droplet is a disk or an annulus and to be given as in \autoref{theo2}, with $\mathcal{F}_{Q}[\mathbb{A}_{r_0,r_1}]$ and $\mathcal{F}_{Q}[\mathbb{D}_{r_1}]$ replaced by $\mathcal{F}_{Q}^{(m)}[\mathbb{A}_{r_0,r_1}]$ and $\mathcal{F}_{Q}^{(m)}[\mathbb{D}_{r_1}]$, respectively, for which an explicit expression is more complicated to derive and somewhat cumbersome.
}
\conj{conj: kernel m}{Consider the same setting as in \autoref{conj: Z_n}. For $k=1,2$, let $z_k=r_*+\frac{\xi_k}{(\gamma_{Q,m}\, n)^{\frac{1}{2m+2}}}\in\mathbb C$, with $\gamma_{Q,m}$ given in \eqref{eq: cst m}, $|\xi_k|\leq \log(n)$. The kernel \eqref{eq: def kernel}, up to a cocycle, admits the following uniform expansion, as $n \to \infty$,
\Eq{}{
&\frac{1}{(\gamma_{Q,m}\, n)^{\frac{1}{m+1}}} K_{n}(z_1,z_2)=\frac{1}{4} K_*^{(m)}(\xi_1,\xi_2) +O(n^{-\frac{1}{2m+2}}),
}
with $K_*^{(m)}$ as in \eqref{eq: kernel m}.
}
One can, consequently, expect the corresponding limiting $1$-point correlation function to be proportional to the function
\Eq{eq: def rho m}{
\rho^{(m)}(\xi):=K_*^{(m)}(\xi,\xi),
}
which, in the same way as \eqref{eq: asymp rho}, can be shown to have the following asymptotic
\Eq{}{
\rho^{(m)}(\xi)\sim \frac{\Re(\xi)^{2m}}{(2m)!},\quad \Re(\xi)\to\infty.
}
Plots of $\rho^{(m)}$, for $m=1,2,3$ are given in \autoref{fig: rho m}. The case $m=1$ corresponds to \autoref{fig: 1pt fct1} as we have $\rho=\rho^{(1)}$.

We believe that the proofs of \autoref{theo1} and \autoref{theo2} can be adapted into proofs of \autoref{conj: Z_n} and \autoref{conj: kernel m}. Namely, one only needs to replace the scale $n^{1/4}$ by $n^{\frac{1}{2m+2}}$ and proceeds with Taylor expansions of the potential up to order $2m+3$ instead of $5$.

\subsection{Example of potential}\label{example potential}
We present here a particular choice of potential verifying the conditions presented in \autoref{pres of model} and for which our main results, consequently, apply. A generalisation of the following example is given in \eqref{eq:cpot}.

Denoting $\abs{z}=r$, for $z\in\mathbb{C}$, the easiest example of potential for which one can observe the emergence of a gap in the droplet, not located at the origin, is
\Eq{eq: example potential}{
Q(z)=Q_{t}(r):=\frac{r^3}{3}-t\frac{r^2}{2}+cr, \quad r,t,c>0.
}
Note that we introduced here a parameter $t$ controlling the formation of the gap. The parameter $c$ does not actually play a role as one can push its dependence in the two parameters $t$ and $n$ via a simple rescaling $x \mapsto x/c$. Thus, without loss of generality we take $c=1$.

\begin{figure}[h!]
        \centering
        \includegraphics[width=0.8\textwidth]{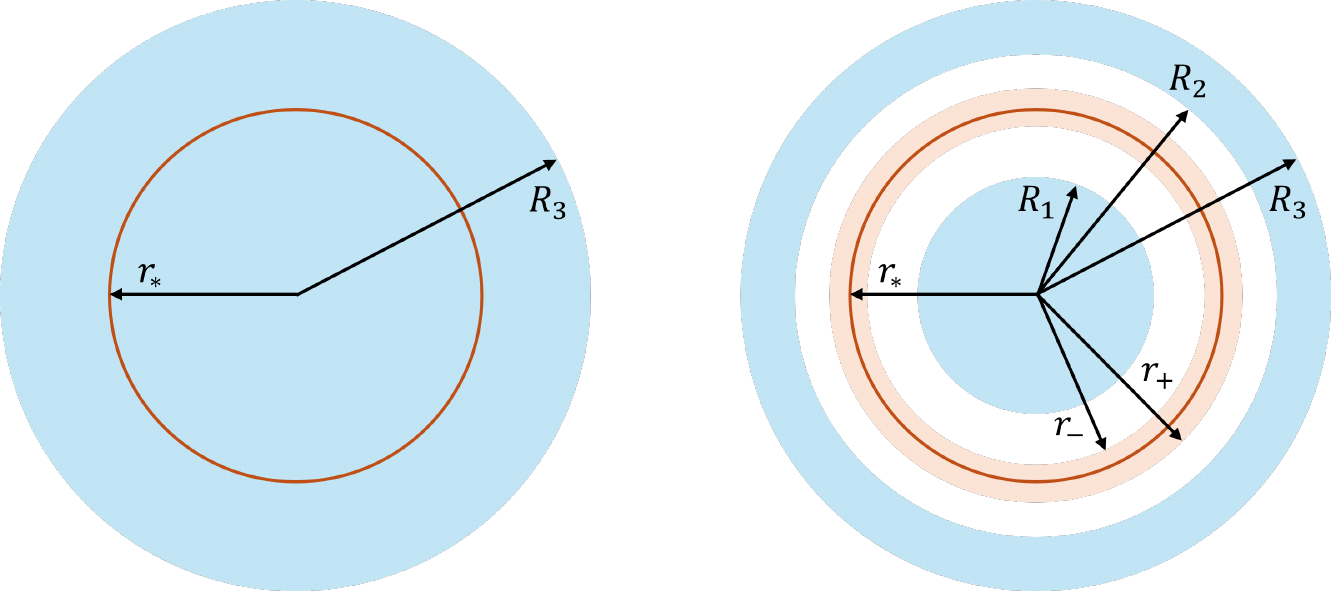}
        \vspace{5pt} 
        \caption{\textit{(Left)} $t=t_{\rm min}$: The droplet (blue) is a disk of radius $R_3$ and both $\Delta Q_{t_{\rm min}}$ and $\partial_r\Delta Q_{t_{\rm min}}$ vanish on the circle of radius $r_*$ (dark orange). \textit{(Right)}  $t>t_{\rm min}$: The droplet (blue) is composed of a disk $\mathbb{D}_{R_1}=\{0\leq \abs{z}\leq R_1 \}$ and an annulus $\mathbb{A}_{R_2,R_3}=\{R_2\leq \abs{z}\leq R_3 \}$. The Laplacian of the potential, $\Delta Q_t$, is strictly negative inside the annulus $\mathbb{A}_{r_-,r_+}=\{r_-\leq \abs{z}\leq r_+ \}$ (orange) and vanishes on the boundaries. }  
        \label{fig: gap apparition}
\end{figure}

Depending on the value of $t$, the support of the equilibrium measure is either a disk, or a disk encircled by a disjoint annulus (cf. \autoref{fig: gap apparition}). There exists a critical value of $t$ for which the annulus starts to separate from the disk. Indeed, the equation
\begin{equation}
\Delta Q_t(r)=Q_t''(r)+\frac{1}{r}Q_t'(r)=0,
\end{equation}
admits the following solutions
\Eq{eq: radii}{
\begin{cases}
    & r_*:=\frac{1}{\sqrt{3}}, \quad t=t_{\rm min}:=\sqrt{3}\\
    & r_\pm:=\frac{t}{3}\pm \sqrt{\frac{1}{3}\left(\frac{t^2}{3}-1\right)},\quad t>t_{\rm min}.
\end{cases}
}

The case $t = t_{\rm min}$ corresponds to the setting described in \autoref{pres of model}, where $r_* = \frac{1}{\sqrt{3}}$ is a double zero of the map $r \mapsto \Delta Q_t(r)$. An illustration of a typical configuration in this case is shown in \autoref{fig: example CG}. For $t < t_{\rm min}$, the potential $Q_{t}$ is strictly subharmonic, and the droplet is a disk. In contrast, for $t > t_{\min}$, $Q_{t}$ is superharmonic in the annulus whose radii are given by $r_{\pm}$ in \eqref{eq: radii}. \autoref{fig: gap apparition} illustrates the configurations $t=t_{\rm min}$ and $t>t_{\rm min}$.

The study of the emergence of a gap by varying the parameter $t$ is the subject of another paper. A multi-scale analysis is required to study this interpolating phenomenon and differs from the one presented in the current work, as the potential $Q_t$ then varies with $n$ via the appropriate $n$-dependent scaling of $t$. 

The rest of the article is organized as follows. In \autoref{preli nota} we introduce notations and lay the groundwork for the large $n$ asymptotic analysis of the norms \eqref{eq:defuj} of the polynomials $p_{j,n}$ which is undertaken in \autoref{sec:asnorm}. In \autoref{sec:theo1} we focus on the asymptotic expansion of the partition function $Z_n$ and give the proof of \autoref{theo1}. In \autoref{local correlations}, we focus on the local correlations of particles in the vicinity of the radius $r_*$ and give, in particular, the proofs of \autoref{theo2}, \autoref{coro2} and \autoref{cor:mls}. Finally, we offer some discussion in \autoref{discussion}.

\section{Preliminaries and notations}\label{preli nota}
To set the stage for the Laplace method used to study the asymptotic expansion of the norms \eqref{eq:defuj}, one can rewrite
\begin{equation}\label{eq: laplace integ}
\frac{\norm{p_{j,n}}^2}{2}=u_j=\int_0^\infty   re^{-nV_\tau(r)}dr
\end{equation}
where we define the auxiliary function $V_\tau$
\begin{equation}
V_\tau(r) := q(r) - 2\tau \log(r), \qquad \tau = \frac{j}{n}.
\end{equation}
The potential $Q$ being subharmonic  implies that the function $r\mapsto rq'(r)$,
for $0\leq\tau\leq1$,
is increasing inside the droplet, i.e. for $r\in[r_0,r_1]$; see \autoref{fig: config}. Therefore, $V_{\tau}$ has a unique minimum $r_\tau$ and the equation $V_\tau'(r_\tau)=0$ can be rewritten
\begin{equation}\label{eq:sadpt}
r_\tau q'(r_\tau) = 2\tau.
\end{equation}
In particular, there exists $\tau_*$ such that  
\begin{equation}\label{eq: tau*}
r_* q'(r_*) = 2\tau_*.
\end{equation}
This $\tau_*$ is special, as the assumption \eqref{eq: vanish laplacian} implies that the first three derivatives of $V_{\tau_*}$ vanish at  $r_*$, i.e.
\begin{equation}
V_{\tau_*}'(r_*) = V_{\tau_*}''(r_*) = V_{\tau_*}'''(r_*) = 0,  
\qquad V_{\tau_*}^{(4)}(r_*) >0.
\end{equation}
A crucial point of the analysis is to control the difference $r_\tau - r_*$. To know how $r_\tau$ varies with $\tau$, one can differentiate the critical point equation \eqref{eq:sadpt} to get the following differential equation
\begin{equation}
    \frac{dr_\tau}{d\tau}=\frac{2}{r_\tau \Delta Q(r_\tau)}.
\end{equation}
Moreover, by the implicit function theorem, one finds
\begin{equation}\label{eq:locexp}
    r_\tau - r_* 
=\left( \frac{12}{r_* \partial^2_r \Delta Q(r_*)} \right)^{1/3}\, (\tau - \tau_*)^{1/3} + o\big( (\tau - \tau_*)^{1/3} \big).
\end{equation}

The difference $r_\tau - r_*$ being an implicit function of the linear parameter $\tau$, the key idea is to reverse the problem and consider instead the difference $\tau-\tau_*$ as an implicit function of the linearly parametrized difference $r_\tau - r_*=\delta_r$, where $\delta_r$ is the new parameter of the problem. This corresponds to the following change of coordinates
\Eq{}{
(\tau, r_\tau)=(\tau_*, r_*)+(p(\delta_r), \delta_r)
}
where the function $p$ is implicitly defined by the relation 
\begin{equation}\label{eq:defp}
p(r_\tau - r_*) = \tau - \tau_*,
\end{equation}
and can be explicitly defined in terms of $q'$ using \eqref{eq:sadpt} and \eqref{eq: tau*},
\begin{equation}\label{eq:defp2}
p(x) := \frac{1}{2}\big[(x+r_*)q'(x+r_*)-r_*q'(r_*)\big].
\end{equation}
The differentiability of $p$ then follows from the differentiability of the potential $q$.

In the forthcoming analysis, we will make use of the following derivatives of $V_\tau$, which we conveniently express in terms of $\Delta Q$,
\Eq{eq: deriv V}{
V_\tau''(r) &= \Delta Q(r) - \frac{1}{r} V_\tau'(r), \\
V_\tau'''(r) &= \partial_r \Delta Q(r) - \frac{1}{r} \Delta Q(r) + \frac{2}{r^2} V_\tau'(r), \\
V_\tau^{(4)}(r) &= \partial_r^2 \Delta Q(r) - \frac{1}{r} \partial_r \Delta Q(r) + \frac{3}{r^2} \Delta Q(r) - \frac{6}{r^3} V_\tau'(r).
}
In a first step, in the large $n$ asymptotic of the norms $\norm{p_{j,n}}$, the above derivatives will be evaluated at $r_\tau$ and thus the $\tau$ dependence on the right-hand side of \eqref{eq: deriv V} will only appear through $r_\tau$, as $ V_\tau'(r_\tau)=0$ by definition. 
In a second step, as we are ultimately interested in the asymptotic of the sum of the norms $\norm{p_{j,n}}$, one needs to determine which indices $j$ will yield the main contribution ($\tau=j/n$). We will see that the main contribution comes from $\tau=\tau_*$ and equivalently from $r_\tau=r_*$---hence the aforementioned change of coordinates.

In this perspective, we give a double expansion of the derivatives of $V_\tau$, depending on how far the argument is from $r_*$ and how far $\tau$ is from $\tau_*$. This will be crucial in our analysis.
\lem{}{
Denoting, $\delta_r=r_\tau-r_*$, the derivatives of $V_\tau$ admit the following expansions, as $\delta_r\to 0$ and $\abs{r-r_*}\to 0$,
\begin{equation}\label{eq:vder}
\begin{aligned}
V_\tau^{(2)}(r) &=  \frac{\gamma_Q}{2}(r-r_*)^2 +\frac{\upsilon_Q}{3!}(r-r_*)^3+\frac{\gamma_Q}{3! r_*}\delta_r^3+ O\left(\abs{\delta_r}^4+\abs{\delta_r}^3\abs{r-r_*}\, +\abs{r-r_*}^4\right), \\
V_\tau^{(3)}(r) &= \gamma_Q \, (r-r_*)  + \frac{\upsilon_Q}{2}(r-r_*)^2+  O\left(\abs{\delta_r}^3+\abs{r-r_*}^3 \right), \\
V_\tau^{(4)}(r) &= \gamma_Q  + \upsilon_Q \, (r-r_*)+ O\left(\abs{\delta_r}^3+\abs{r-r_*}^2 \right), \\
V_\tau^{(5)}(r) &=\upsilon_Q+ \, O\left(\abs{\delta_r}^3+\abs{r-r_*} \right),
\end{aligned}
\end{equation}
with 
\Eq{eq: upsilonQ}{\gamma_Q:=\partial_r^2\Delta Q(r_*),
\qquad
\upsilon_Q:=\partial_r^3\Delta Q(r_*)-\frac{\gamma_Q}{r_*}.} 
}
\begin{proof}
The first step is to proceed with a Taylor expansion around $r=r_*$ of the successive derivatives of $V_\tau$, i.e.
\Eq{eq: taylor deriv}{
V_\tau^{(k)}(r)=V_\tau^{(k)}(r_*)+V_\tau^{(k+1)}(r_*)(r-r_*)+\frac{V_\tau^{(k+2)}(r_*)}{2}(r-r_*)^2+...
}
We recall that, by assumption, $V_\tau$ is at least six times continuously differentiable (see \autoref{pres of model}), which ensures that a first-order Taylor expansion of $V_\tau^{(5)}$ is well defined.
Then, we use \eqref{eq: deriv V} and recall that, by assumption, 
\Eq{eq: assum}{
& \Delta Q(r_*)=q''(r_*)+\frac{q'(r_*)}{r_*}=0,\quad
\partial_r\Delta Q(r_*) = q'''(r_*) + \frac{q''(r_*)}{r} - \frac{q'(r_*)}{r^2}=0
}
with the first non-vanishing derivative given by
\Eq{}{
\gamma_Q=\partial^2_r \Delta Q(r_*)
&= q^{(4)}(r_*) + \frac{q^{(3)}(r_*)}{r_*}  - \frac{2q''(r_*)}{r_*^2} + \frac{2q'(r_*)}{r_*^3}>0,
}
which is equal to \eqref{eq:vcon}, using \eqref{eq: assum}.
As the successive derivatives of $V_\tau$ are all expressed in terms of $V_\tau'$, from \eqref{eq: deriv V}, one can further use the expansion \eqref{eq:locexp} (or, equivalently Taylor expand the left-hand side of \eqref{eq:defp} using \eqref{eq:defp2}) to get
\Eq{eq: expansion v prime}{
 V_\tau'(r_*)=-2\frac{\tau-\tau_*}{r_*}=-\frac{\gamma_Q}{3!} \delta_r^3+O(\delta_r^4).
}
Combining \eqref{eq: expansion v prime} and \eqref{eq: taylor deriv} with the help of \eqref{eq: deriv V} yields the claim.

\end{proof}

\nota{nota: leq geq}{To avoid excessive use of unspecified constants, we write “$A_n \lesssim B_n$” to
denote that $A_n \leq C B_n$ for all large enough $n$, where $C$ is some constant independent of $n$. 
}

\section{Asymptotic expansion of the norms $\norm{p_{j,n}}$}\label{sec:asnorm}
In this section, we study the large $n$ asymptotics of the norms $\norm{p_{j,n}}^2=2 u_j$ \eqref{eq: norms} which appear both in the kernel \eqref{eq: def kernel} and in the expression of the partition function \eqref{eq:defuj}. 

While we treat $\tau$ as a continuous parameter, we remind that $\tau=j/n$ and thus we define 
\Eq{}{
j_*:= n\tau_*,
}
which is not necessarily an integer. As we will consider two different regimes, depending on the size of the difference $\tau-\tau_*$ and ultimately of the difference $j-j_*$, we introduce the following scalings
\begin{equation}\label{eq:defdn}
    [-\ell_n^{-},\ell_n^{+}]:= n \,p\!\left( \left[-\frac{\xi_n}{n^{1/4}},\frac{\xi_n}{n^{1/4}}\right] \right),
\end{equation}
where $\xi_n$ has to be carefully chosen such that
\Eq{eq: prop xi}{
\frac{n^{1/4}}{\xi_n^5}\to 0,\qquad \frac{\xi_n^3}{n^{1/4}}\to 0, \qquad n\to\infty.
}
We then choose $\xi_n$ to be
\Eq{eq: def xi}{
\xi_n:=n^{1/15}.
}
The scalings $\ell_n^{\pm}$ are defined implicitly by the map $p$ given in \eqref{eq:defp2}.

The first regime is for $j$ close to $j_*$ ($\tau$ close to $\tau_*$) and the corresponding asymptotics of the norms $u_j$ is given in the following proposition.
\prop{prop:norm}{
For $j\in  [j_*-\ell_n^{-},j_*+\ell_n^{+}]$, the norms $u_j$ \eqref{eq:defuj} admit the following uniform asymptotic expansion, as $n \to \infty$, 
\Eq{eq:ujasymp1}{
u_j &= \frac{e^{-n V_\tau(r_\tau)} r_\tau}{n^{1/4}}
\left( f_1(x) + \frac{1}{n^{1/4}} f_2(x) + O\!\left(\frac{1}{n^{1/2}}\right) \right) \\
&= \sqrt{\frac{2\pi}{n V_\tau''(r_\tau)}} \, r_\tau \, e^{-n V_\tau(r_\tau)}\sqrt{\frac{\gamma_Q}{4\pi}} \, |x| \, f_1(x)
\left(1 + \frac{1}{n^{1/4}} \left[\frac{f_2(x)}{f_1(x)}
+
\frac{ \left(\upsilon_Q+\frac{\gamma_Q}{r_*}\right)}{6\gamma_Q}\,x\right]+ O\!\left( \frac{1}{n^{1/2}} \right) \right),
}
where $x = n^{1/4} (r_\tau-r_*)$ (which depends on $j$ through $\tau$) and 
\Eq{eq:f1}{
f_1(x) :=&\, \int_{-\infty}^{\infty} \exp\!\left[-\gamma_Q\left(
   x^2\frac{ y^2}{4} \, 
 + x\frac{   y^3}{3!} \,
+ \frac{ y^4}{4!} \,\right)
\right] \, dy,\\
f_2(x)
:=&\,
\int_{-\infty}^{\infty}
\left[\frac{y}{r_*}-\left(\upsilon_Q+\frac{\gamma_Q}{r_*}\right)\frac{x^3}{3!}\frac{ y^2}{2}
-\upsilon_Q\left(\frac{x^2}{2}\frac{y^3}{3!}
+x\frac{y^4}{4!}
+\frac{y^5}{5!}\right)\right]\\
&\qquad\times\exp\!\left[-\gamma_Q\left(
   x^2\frac{ y^2}{4} \, 
 + x\frac{   y^3}{3!} \,
+ \frac{ y^4}{4!} \,\right)
\right] \, dy,
}
with $\gamma_Q = \partial^2_r \Delta Q(r_*)$ and $\upsilon_Q=\partial_r^3\Delta Q(r_*)-\frac{\gamma_Q}{r_*}$.
}

\begin{proof}
Let us proceed with Laplace's method on the integral 
\Eq{}{
u_j=\int_0^\infty   re^{-nV_\tau(r)}dr.
}
Due to the large $n$ parameter, the main contribution of the above integral will come from a small interval $[r_\tau-\beta_n, r_\tau+\beta_n]$ centred around the minimum of the function $V_\tau$, which is $r_\tau$. The cut-off scaling $\beta_n$ will be shortly determined. On this interval, a good approximation of $V_\tau$ is given by a truncated Taylor expansion. As pointed out in \autoref{preli nota}, the function $V_\tau$ vanishes up to third order for $(\tau,r_\tau)=(\tau_*,r_*)$. Thus, one needs to Taylor expand $V_\tau$ around $r=r_\tau$ up to fifth order, to control the error. 

Let us start by decomposing the integral $u_j$ as follows
\begin{equation}
u_j = I_1 + I_2 + I_3,
\end{equation}
where
\begin{equation}
I_1 = \int_{0}^{r_\tau - \beta_n} r e^{-n V_\tau(r)} \, dr, \quad
I_2 = \int_{r_\tau - \beta_n}^{r_\tau + \beta_n} r e^{-n V_\tau(r)} \, dr, \quad
I_3 = \int_{r_\tau + \beta_n}^{\infty} r e^{-n V_\tau(r)} \, dr.
\end{equation}
We will show that the tail contributions $I_1$ and $I_3$ are negligible compared to $I_2$. 

Regarding the cut-off scaling $\beta_n$, it should be slightly larger than the scales of $r-r_\tau$ that keeps the difference $n(V_\tau(r)-V_\tau(r_\tau))$ exactly of order $1$, independently of $n$ and $\tau$ (in the range considered). More precisely, if one proceeds with the change of variables $r\mapsto r_\tau+r/c_n$, an appropriate interpolating scale $c_n$ is given by
\begin{equation} \label{eq:defcn}
c_n := n^{1/2} \abs{r_\tau-r_*}+ n^{1/4},
\end{equation}
as can be seen by Taylor expanding $V_\tau$ around $r_\tau$. From there, one finds that a suitable cut-off scale $\beta_n$ is
\begin{equation} \label{eq:defbetan}
\beta_n := \frac{\log(n)}{c_n}.
\end{equation}

\textbf{Step 1: Control of the tails.} Let us first deal with the contributions coming from $I_1$ and $I_3$. Since $V_\tau$ has a unique minimum at $r_\tau$, it is increasing above that point. Thus, to estimate the tail contribution for $r \geq r_\tau + \beta_n$ one can use a Taylor expansion of $V_\tau$ to get
\begin{equation}
V_\tau(r) - V_\tau(r_\tau) 
\ \geq\ V_\tau(r_\tau + \beta_n) - V_\tau(r_\tau) 
\ \geq\  \frac{V_\tau^{(2)}(r_\tau)}{2} \beta_n^2 
 + \frac{V_\tau^{(3)}(r_\tau)}{3!} \beta_n^3 
 + \frac{V_\tau^{(4)}(r_\tau)}{4!} \beta_n^4+O(\beta_n^5).
\end{equation}
Since the coefficients still depend on $n$ through $\tau$ and that $V_\tau^{(2)}(r_\tau)$ and $V_\tau^{(3)}(r_\tau)$ vanish at $\tau = \tau_*$, we may use the definition of $\beta_n$ in \eqref{eq:defbetan} to take a looser lower bound:
\Eq{eq: bound diff}{
V_\tau(r) - V_\tau(r_\tau)\geq C \frac{\log(n)^2}{n},
}
where $C$ is some positive constant independent of $n$. 

Then, using \eqref{eq: growth pot}, we choose a sufficiently large $M > 0$, independent of $n$, such that
\Eq{}{
\int_{0}^{\infty} \left( r^{2} e^{-q(r)} \right)^{M} r \, dr < +\infty,
}
and combine this with \eqref{eq: bound diff} to obtain the upper bound
\begin{equation}
I_3 \leq e^{-n V_\tau(r_\tau)}e^{- C(n-M)\frac{\log(n)^2}{n}}\int_{r_\tau + \beta_n}^{\infty} r e^{-M (V_\tau(r)-V_\tau(r_\tau))} \, dr \lesssim e^{-n V_\tau(r_\tau)} \, e^{- C\log(n)^2},
\end{equation}
which is exponentially small in $n$. The notation $\lesssim$ is explained in \autoref{nota: leq geq}. A similar argument for $r \leq r_\tau - \beta_n$ yields
\begin{equation}
I_1\lesssim e^{-n V_\tau(r_\tau)} \, e^{- C\log(n)^2}.
\end{equation}
Thus both $I_1$ and $I_3$ are negligible in comparison with $I_2$.
\medskip
\noindent

\textbf{Step 2: Local expansion near $r_\tau$.}  
We now focus on $I_2$. For $|r - r_\tau| \leq \beta_n$, the Taylor expansion of $V_\tau$ about $r_\tau$ up to fifth order with mean value remainder gives
\begin{equation}
\begin{aligned}
V_\tau(r) - V_\tau(r_\tau) 
= \frac{V_\tau^{(2)}(r_\tau)}{2} (r - r_\tau)^2 
 &+ \frac{V_\tau^{(3)}(r_\tau)}{3!} (r - r_\tau)^3 
 + \frac{V_\tau^{(4)}(r_\tau)}{4!} (r - r_\tau)^4 \\
&\quad + \frac{V_\tau^{(5)}(r_\tau)}{5!} (r - r_\tau)^5+O\left((r - r_\tau)^6\right),
\end{aligned}
\end{equation}
Substituting into $I_2$, we obtain
\begin{equation}
\begin{aligned}
I_2 &= e^{-n V_\tau(r_\tau)}
\int_{r_\tau - \beta_n}^{r_\tau + \beta_n} 
 r \, \exp\!\left[
 -n\left(
 \frac{V_\tau^{(2)}(r_\tau)}{2} (r - r_\tau)^2
 + \frac{V_\tau^{(3)}(r_\tau)}{3!} (r - r_\tau)^3 \right.\right.\\
&\quad\left.\left.
 + \frac{V_\tau^{(4)}(r_\tau)}{4!} (r - r_\tau)^4
 + \frac{V_\tau^{(5)}(r_\tau)}{5!} (r - r_\tau)^5+O\left((r - r_\tau)^6\right)
 \right)
 \right] dr .
\end{aligned}
\end{equation}
\medskip
\noindent
Proceeding with the change of variable $r= r_\tau+ \frac{y}{n^{1/4}}$ yields
\Eq{eq: I2 intermediate}{
I_2 =&\, \frac{e^{-n V_\tau(r_\tau)}}{n^{1/4}} 
\int_{-n^{1/4}\beta_n}^{n^{1/4}\beta_n} 
\left(r_\tau + \frac{y}{n^{1/4}}\right)
\exp\!\Bigg[
 -n\Bigg(
 \frac{V_\tau^{(2)}(r_\tau)}{2} \left(\frac{y}{n^{1/4}}\right)^2
 + \frac{V_\tau^{(3)}(r_\tau)}{3!} \left(\frac{y}{n^{1/4}}\right)^3 \\
&
 + \frac{V_\tau^{(4)}(r_\tau)}{4!} \left(\frac{y}{n^{1/4}}\right)^4
 + \frac{V_\tau^{(5)}(r_\tau)}{5!} \left(\frac{y}{n^{1/4}}\right)^5+O\left( \left(\frac{y}{n^{1/4}}\right)^6\right)
 \Bigg)
 \Bigg] dy .
}
Thus, using the expansion \eqref{eq:vder} and rescaling the difference $r_\tau-r_*$ as follows
\begin{equation}\label{eq:scaling}
\delta_r =r_\tau-r_*= \frac{x}{n^{1/4}},
\end{equation}
we find 
\Eq{eq: phase expansion}{
&n\Bigg(
 \frac{V_\tau^{(2)}(r_\tau)}{2} \left(\frac{y}{n^{1/4}}\right)^2
 + \frac{V_\tau^{(3)}(r_\tau)}{3!} \left(\frac{y}{n^{1/4}}\right)^3 
 + \frac{V_\tau^{(4)}(r_\tau)}{4!} \left(\frac{y}{n^{1/4}}\right)^4
 + \frac{V_\tau^{(5)}(r_\tau)}{5!} \left(\frac{y}{n^{1/4}}\right)^5 \Bigg)\\
=&\, \gamma_Q\left(\frac{ x^2}{2} \, \frac{y^2}{2}
 + x \, \frac{y^3}{3!}
 + \,\frac{y^4}{4!}\right)+\frac{1}{n^{1/4}}\left(\left(\upsilon_Q+\frac{\gamma_Q}{r_*}\right)\frac{x^3}{3!}\frac{y^2}{2}
+\upsilon_Q\left[\frac{x^2}{2}\frac{y^3}{3!}
+x\frac{y^4}{4!}
+\frac{y^5}{5!}\right]\right)\\
&\,+O(n^{-1/2} ),
}
where the right-hand side of the above equation can be written in the more compact form
\Eq{}{
\gamma_Q x^4 \mathbf{P}'\left(\frac{y}{x}\right)+\frac{1}{n^{1/4}}\left(\left(\upsilon_Q+\frac{\gamma_Q}{r_*}\right)\frac{x^3}{3!}\frac{y^2}{2}
+\upsilon_Q x^5 \mathbf{P}\left(\frac{y}{x}\right)\right)+O(n^{-1/2} )
}
with 
\Eq{eq: def bf P}{
\mathbf{P}\left(x\right):=\frac{x^3}{2(3!)}+\frac{x^4}{4!}+\frac{x^5}{5!},\quad \mathbf{P}'(x)=P(x).
}
Plugging back \eqref{eq: phase expansion} in \eqref{eq: I2 intermediate} and expanding the exponential as follows
\Eq{}{
&\exp[-\gamma_Q x^4 \mathbf{P}'\left(\frac{y}{x}\right)-\frac{1}{n^{1/4}}\left(\left(\upsilon_Q+\frac{\gamma_Q}{r_*}\right)\frac{x^3}{3!}\frac{y^2}{2}
+\upsilon_Q x^5 \mathbf{P}\left(\frac{y}{x}\right)\right)+O(n^{-1/2} )]\\
=&\,\exp[-\gamma_Q x^4 \mathbf{P}'\left(\frac{y}{x}\right)]\left(1-\frac{1}{n^{1/4}}\left(\left(\upsilon_Q+\frac{\gamma_Q}{r_*}\right)\frac{x^3}{3!}\frac{y^2}{2}
+\upsilon_Q x^5 \mathbf{P}\left(\frac{y}{x}\right)\right)+O(n^{-1/2} )\right)
}
with the help of the expansion 
\Eq{}{
\frac1{r_\tau}=\frac1{r_*}+O(n^{-1/4}),
}
yields
\Eq{}{
I_2 = \frac{e^{-n V_\tau(r_\tau)}r_\tau}{n^{1/4}} 
\int_{-\frac{\log(n)}{1+\abs{x}}}^{\frac{\log(n)}{1+\abs{x}}} 
&\left(1 + \frac{y}{r_* n^{1/4}}\right)
\exp[-\gamma_Q x^4 \mathbf{P}'\left(\frac{y}{x}\right)]\\
&\times\left(1-\frac{1}{n^{1/4}}\left(\left(\upsilon_Q+\frac{\gamma_Q}{r_*}\right)\frac{x^3}{3!}\frac{y^2}{2}
+\upsilon_Q x^5 \mathbf{P}\left(\frac{y}{x}\right)\right)+O(n^{-1/2} )\right)dy .
}
Finally, given that the integrand in $I_2$ decays exponentially as $|y| \to \infty$, we extend the integration range from $[-\log(n), \log(n)]$ to $\mathbb{R}$, incurring only an exponentially small error in $n$ and we get the first line of \eqref{eq:ujasymp1}
\Eq{eq: uj intermediate}{
u_j &= \frac{e^{-n V_\tau(r_\tau)} r_\tau}{n^{1/4}}
\left( f_1(x) + \frac{1}{n^{1/4}} f_2(x) + O\!\left(\frac{1}{n^{1/2}}\right) \right)\\
&= \sqrt{\frac{2\pi}{n V_\tau''(r_\tau)}} \, r_\tau \, e^{-n V_\tau(r_\tau)}
\left( n^{1/4} \sqrt{\frac{V_\tau''(r_\tau)}{2\pi}}\left[f_1(x) + \frac{1}{n^{1/4}} f_2(x) + O\!\left(\frac{1}{n^{1/2}}\right)\right] \right),
}
where $f_1$ and $f_2$ are given in \eqref{eq:f1}. 

To derive the second line of \eqref{eq:ujasymp1}, one uses the first equation in \eqref{eq:vder} along with \eqref{eq:scaling} to find
\begin{equation}
n^{1/4}\sqrt{\frac{V_\tau''(r_\tau)}{2\pi}}
=
\sqrt{\frac{\gamma_Q}{4\pi}}\,|x|
+\frac{\left(\upsilon_Q+\frac{\gamma_Q}{r_*}\right)}{12}\sqrt{\frac{1}{\pi\gamma_Q}}\,x|x|\,\frac{1}{n^{1/4}}
+O(n^{-1/2}).
\end{equation}
Plugging the above expansion back in the second line of \eqref{eq: uj intermediate} and collecting the terms finishes the proof.
\end{proof}

The second regime is  for $j$ sufficiently far from $j_*$ ($\tau$ far from $\tau_*$). In this case, $V_\tau''(r_\tau)$ does not vanish. However, one needs to distinguish two cases: the case where the droplet is an annulus ($r_0>0$) and the case where it is a disk ($r_0=0$).

\subsection{Annulus case}
\prop{lem:norm}{
   Let $r_0>0$ (cf. \autoref{fig: config} left) and $j\notin  [j_*-\ell_n^{-},j_*+\ell_n^{+}]$, then as $n\rightarrow\infty$ we have 
    \begin{equation}
       u_j=\sqrt{\frac{2\pi}{n V_\tau''(r_\tau)}}r_\tau e^{-n V_\tau(r_\tau)}\left(1+\frac{1}{n} \mathcal{B}(r_\tau)+O(n^{-2})\right)
    \end{equation}
    where 
\Eq{eq: def B}{
\mathcal{B}(r):=- \frac{1}{8} \frac{\partial_r^2 \Delta Q(r)}{\left(\Delta Q(r)\right)^2}
- \frac{19}{24 r} \frac{\partial_r \Delta Q(r)}{\left(\Delta Q(r)\right)^2}
+ \frac{5}{24} \frac{\left(\partial_r \Delta Q(r)\right)^2}{\left(\Delta Q(r)\right)^3}
+ \frac{1}{3r^2} \frac{1}{\Delta Q(r)}.
}
}
\begin{proof}

As $ V_\tau''(r_\tau)$ is not vanishing in this regime, the proof follows from the standard Laplace's method, see e.g. \cite[p.37]{Temme2014}. The same computations (up to a different cutoff scaling, on which the result does not depend) are carried out in \cite[proof of Lemma 3.2]{BKS23}. We note that, in that reference, $\Delta$ is defined as one quarter of the usual Laplacian.
\end{proof}
\subsection{Disk case}
When considering a disk droplet (i.e., $r_0=0$), the saddle point $r_\tau$ comes close to  the boundary of the integration domain. This situation requires a refined local analysis, which is carried out in \cite[\S 3.1]{BKS23}, for indices going up to order $O(n^{\epsilon})$ for some $0<\epsilon<1/5$. With this choice, the following large $n$ asymptotics for $u_j$ holds.

\prop{lem: uj asym6}{(\cite[Lemma 3.1, specialized to $h_j=2u_j$]{BKS23})  
If $r_0=0$ (cf. \autoref{fig: config} right) and $j \leq n^{\epsilon}$ with $0<\epsilon<1/5$, then as $n\to\infty$ the integral $u_j$ defined in \eqref{eq:defuj} satisfies  
\begin{equation}\label{eq: uj asym6}
     2 u_j = -2n q(0) - (j+1)\log\!\big(n q''(0)\big) + \log(j!) 
     + {\rm O}\!\left(j^{3/2}\frac{\log (n)^3}{\sqrt{n}}\right).
\end{equation}
}

\section{Asymptotic expansion of $Z_n$}\label{sec:theo1} 
Using the asymptotic of the norms $u_j$ we can now compute the large $n$ asymptotic of the partition function $Z_n$. 

\begin{proof}[Proof of \autoref{theo1}]
Let us first split the sum as follows
\begin{equation}
   \log(\frac{Z_n}{(2\pi)^n })= \sum_{j=0}^{n-1} \log(u_j) = \sum_{j=0}^{\lfloor j_* - \ell_n^{-} \rfloor - 1} \log(u_j)
    + \sum_{j=\lfloor j_* + \ell_n^{+} \rfloor + 1}^{n-1} \log(u_j)+ \sum_{j=\lfloor j_* - \ell_n^{-} \rfloor}^{\lfloor j_* + \ell_n^{+} \rfloor} \log(u_j),
\end{equation}
where $\lfloor \cdot \rfloor$ is the floor function and $\ell_n^\pm$ as in \eqref{eq:defdn}. Since the asymptotic behavior of the norms $u_j$ depends on whether the droplet is a disk or an annulus, the same distinction applies here.
\subsection{Annulus case}
We remind that $\tau=j/n$ and that $x = n^{1/4} (r_\tau-r_*)$ (which depends on $j$ through $\tau$). From \autoref{prop:norm}, we have for  
$j \in [\,\lfloor j_* - \ell_n^{-} \rfloor, \, \lfloor j_* + \ell_n^{+} \rfloor\,]$,  
\begin{equation}\label{eq:sum2}
    \log(u_j)
    = -n V_\tau(r_\tau)
      - \frac12 \log(n)
      + \frac12 \log\!\left( \frac{2\pi r_\tau^2}{V_\tau''(r_\tau)} \right)
      + \log(h_1(x))
      +\frac{1}{n^{1/4}}h_2(x)
      + O(n^{-1/2}),
\end{equation}
where we define
\Eq{eq:defh1}{
h_1(x):=\sqrt{\frac{\gamma_Q}{4\pi}} \, |x| \, f_1(x), \qquad h_2(x):=\frac{f_2(x)}{f_1(x)}
+
\frac{ \left(\upsilon_Q+\frac{\gamma_Q}{r_*}\right)}{6\gamma_Q}\,x
}
with $f_1$ and $f_2$ given in \eqref{eq:f1}. 

On the other hand, in the annulus case $(r_0 > 0)$, it follows from \autoref{lem:norm} that for
$j \notin [\,\lfloor j_* - \ell_n^{-} \rfloor,\, \lfloor j_* + \ell_n^{+} \rfloor\,]$,
\begin{equation}\label{eq:sum1}
    \log(u_j)
    = -n V_\tau(r_\tau)
      - \frac12 \log(n)
      + \frac12 \log\!\left( \frac{2\pi r_\tau^2}{V_\tau''(r_\tau)} \right)
      +\frac{1}{n} \mathcal{B}(r_\tau)+ O(n^{-2}).
\end{equation}
Combining \eqref{eq:sum2} and \eqref{eq:sum1} yields  
\begin{equation}
    \sum_{j=0}^{n-1} \log(u_j)
    = T_1 + T_2 +T_3+T_4+ O\left(\frac{\xi_n^3}{n^{1/4}}\right), 
\end{equation}
where the error term $O(\xi_n^3 n^{-1/4})=o(1)$ follows from the definition \eqref{eq:defdn} combined with \eqref{eq:defp} and \eqref{eq:locexp}, and
\Eq{}{
    T_1
    &:= \sum_{j=0}^{n-1}
      \left[
        -n V_\tau(r_\tau)
        - \frac12 \log(n)
        + \frac12 \log\!\left( \frac{2\pi r_\tau^2}{V_\tau''(r_\tau)} \right)
      \right], \qquad
      T_2 := \sum_{j=\lfloor j_* - \ell_n^{-} \rfloor}^{\lfloor j_* + \ell_n^{+} \rfloor} \log(h_1(x)),\\
    T_3&:=\sum_{j=\lfloor j_* - \ell_n^{-} \rfloor}^{\lfloor j_* + \ell_n^{+} \rfloor} \frac{1}{n^{1/4}}h_2(x) ,\qquad
    T_4:=\sum_{j \notin [\,\lfloor j_* - \ell_n^{-} \rfloor,\, \lfloor j_* + \ell_n^{+} \rfloor\,]} \frac{1}{n}\mathcal{B}(r_\tau).
}
The large $n$ asymptotics of the different $T_k$, $k={1,\cdots,4}$, are then given in the following lemmas.
\lem{ref:sumT1}{ As $n\rightarrow\infty$, we have 
\begin{equation}\label{eq:sumt1}\begin{aligned}
    T_1
    =&\, \sum_{j=0}^{n-1}
      \left[
        -n V_\tau(r_\tau)
        - \frac12 \log(n)
        + \frac12 \log\!\left( \frac{2\pi r_\tau^2}{V_\tau''(r_\tau)} \right)
      \right],\\
      =&\,-n^2 I_Q[\mu_Q]-\frac{n\log(n)}{2}+n\left[\frac{\log (\pi/2)}{2}-\frac{1}{2}E_Q[\mu_Q]\right]\\
      &\,-\frac{1}{3}\log\!\left(\frac{r_1}{r_0}\right)
+ \frac{1}{4}\log\!\left(\frac{\Delta Q(r_1)}{\Delta Q(r_0)}\right)+O(n^{-1}).
\end{aligned}\end{equation}}

\begin{proof} The formula \eqref{eq:sumt1} follows from a straightforward application of the Euler-Maclaurin formula given by \autoref{prop: EM formula}; see \cite[Lemma 2.3 and Lemma 2.4]{BKS23} or \cite[Proposition 4.2]{AFLS25}.

\end{proof}
\lem{ref:sumT2}{
As $n \to \infty$, we have 
\begin{equation}
    T_2
    = \sum_{j=\lfloor j_* - \ell_n^{-} \rfloor}^{\lfloor j_* + \ell_n^{+} \rfloor} \log( h_1\!\left(x\left(j/n\right)\right))
    = n^{1/4} \frac{\varkappa_Q}{2}\int_{-\infty}^{\infty}x^2\log(\frac{x^2}{\sqrt{4\pi}} \mathcal{P}(x^4))dx -\frac{7 r_*}{12}\frac{n^{1/4}}{\xi_n}+ o(1),
\end{equation}
where the function $h_1$ is defined in~\eqref{eq:defh1} and $\varkappa_Q$, $\mathcal{P}$ and $\xi_n$ are given, respectively, in \eqref{eq:vcon}, \eqref{eq: pearcey} and \eqref{eq: def xi}.
}

\begin{proof}
Before applying Euler-Maclaurin formula, one needs to determine the large argument asymptotics of the function $x\mapsto \log(h_1(x))$. By proceeding with the change of variable $y\mapsto xy$ in the integral expression of $f_1$ \eqref{eq:f1} one gets
\Eq{}{
\log(h_1\left(x\gamma_Q^{-1/4}\right)))=\log(\frac{x^2}{\sqrt{4\pi}} \mathcal{P}(x^4)).
}
Then, from the standard Laplace method, one finds that $\mathcal{P}(\theta)$, defined in \eqref{eq: pearcey}, has the large $\theta$ asymptotic behaviour
\begin{equation}
    \mathcal{P}(\theta)=\sqrt{\frac{4\pi}{\theta}}\left(1+\frac{7}{6\theta}+O\left(\frac{1}{\theta^2}\right)\right),\qquad \theta\rightarrow\infty. 
\end{equation}
Thus, it follows that
\begin{equation}\label{eq:asy}
    \log(\frac{x^2}{\sqrt{4\pi}} \mathcal{P}(x^4))=\frac{7}{6x^4}+O\left(\frac{1}{x^8}\right),\qquad x\to \pm\infty.
\end{equation}

We are now ready to apply the Euler-Maclaurin formula (cf. \autoref{prop: EM formula}) on $T_2$. Reminding that the variable $x$ is implicitly a function of $\tau=j/n$, we note $x=x(\tau)$ and we get
\Eq{eq: EM T2}{
    \sum_{j = \lfloor j_* - \ell_n^{-} \rfloor}^{\lfloor j_* + \ell_n^{+} \rfloor}
        \log( h_1\!\left(x\left(j/n\right)\right))
    =&\, n \int_{\tau_n^{-}}^{\tau_n^{+}}
        \log( h_1(x(\tau))) \, d\tau\\ &+ \frac{\log( h_1(x(\tau_n^{-})))}{2}+  \frac{\log( h_1(x(\tau_n^{+})))}{2}+O(n^{-1}),
}
with the sequence
\Eq{}{
\tau_n^{\pm}:=\lfloor j_* \pm \ell_n^\pm \rfloor/n.
}
Starting with the boundary terms in the right-hand side of \eqref{eq: EM T2}, one uses that
\Eq{}{
x(\tau_n^\pm)\sim \pm\xi_n,\qquad n\to\infty,
}
due to \eqref{eq:defdn}, along with the asymptotic \eqref{eq:asy}, to find that they are negligible as
\Eq{}{
\log( h_1(x(\tau_n^{\pm})))=O\left(\xi_n^{-4}\right),\qquad n\to \infty.
}
Turning to the first term on the right-hand side of \eqref{eq: EM T2}, in preparation of a subsequent change of variables, one can replace $\tau_n^{\pm}$ by
\Eq{}{
\tilde{\tau}_n^{\pm}:= (j_* \pm \ell_n^\pm) /n,
}
incurring a negligible error of order  
\Eq{}{
O\!\left(n|\tau_n^\pm-\widetilde\tau_n^\pm|\sup_{|x|\ge \xi_n}|\log( h_1(x))|\right)
=O\left(\xi_n^{-4}\right),\qquad n\to \infty.
}
Thus, proceeding with the change of variable  
\Eq{}{
    \tau= \tau_* +p\!\left(\frac{x}{n^{1/4}}\right),
}
and reminding the definition of $\ell_n^\pm$ \eqref{eq:defdn}, the integral becomes
\begin{equation}\label{hint}
     n \int_{\tilde{\tau}_n^{-}}^{\tilde{\tau}_n^{+}}
        \log( h_1(x(\tau))) \, d\tau
    = n^{3/4}  \int_{-\xi_n}^{\xi_n}
        \log( h_1(x)) \, p'\!\left(\frac{x}{n^{1/4}}\right) dx.
\end{equation}
The explicit expression of $p'$ follows from \eqref{eq:defp2},  
\begin{equation}
    p'\!\left(\frac{x}{n^{1/4}}\right)
    = \frac{\left(r_* + \frac{x}{n^{1/4}}\right)
             q''\!\left(r_* + \frac{x}{n^{1/4}}\right)
           + q'\!\left(r_* + \frac{x}{n^{1/4}}\right)}{2},
\end{equation}
and an asymptotic expansion, as $n\to\infty$, is given by  
\Eq{eq: asymp p'}{
 p'\!\left(\frac{x}{n^{1/4}}\right)
 &= \frac{r_*  q^{(4)}(r_*) + 3  q^{(3)}(r_*)}{4\, n^{1/2}}x^2
    + \frac{r_*  q^{(5)}(r_*) + 4 q^{(4)}(r_*)}{12\, n^{3/4}} x^3
    + O(n^{-1}) \\
 &= \frac{r_*\gamma_Q }{4\,n^{1/2}}x^2+ \frac{r_*\upsilon_Q + 4 \gamma_Q}{12\, n^{3/4}} x^3
    + O(n^{-1}),
}
where we have used \eqref{eq: assum}. Proceeding with the change of variables $y\mapsto xy$ in the integral expression of $f_1$ \eqref{eq:f1} and then $x\mapsto x/\gamma_Q^{1/4}$ in \eqref{hint} yields
\Eq{eq: T2 integ}{
T_2=n^{1/4} \frac{\varkappa_Q}{2} \int_{-\xi_n\gamma_Q^{1/4}}^{\xi_n\gamma_Q^{1/4}}x^2\log(\frac{x^2}{\sqrt{4\pi}} \mathcal{P}(x^4))dx
      + O\left(\xi_n^{-4}\right).
}
Note that the term of order $O(n^{-3/4})$, in the expansion \eqref{eq: asymp p'}, is associated with an odd power of $x$ and yields a vanishing contribution, hence the error term $O(\xi_n^{-4})$ in the above equation.

Finally, as the integrand in \eqref{eq: T2 integ} is independent of $n$ and integrable on $\mathbb{R}$ due to \eqref{eq:asy}, the bounds $\pm \xi_n\gamma_Q^{1/4}$ can be extended to $\pm \infty$ while controlling the error by integrating the asymptotic expansion of the integrand \eqref{eq:asy}. Explicitly, 
\Eq{}{
T_2&=n^{1/4} \frac{\varkappa_Q}{2} \left(\int_{-\infty}^{\infty}-2\int_{\xi_n\gamma_Q^{1/4}}^\infty\right)x^2\log(\frac{x^2}{\sqrt{4\pi}} \mathcal{P}(x^4))dx + O\left(\xi_n^{-4}\right)\\
&=n^{1/4} \frac{\varkappa_Q}{2}\int_{-\infty}^{\infty}x^2\log(\frac{x^2}{\sqrt{4\pi}} \mathcal{P}(x^4))dx -\frac{7}{12}r_*\frac{n^{1/4}}{\xi_n}+ O\left(\frac{n^{1/4}}{\xi_n^5}\right),
}
where the error verifies
\Eq{}{
O\left(\frac{n^{1/4}}{\xi_n^5}\right)=o(1),\qquad n\to\infty,
}
due to the suitable choice \eqref{eq: def xi}. This finishes the proof.
\end{proof}

\lem{ref:sumT3}{
As $n \to \infty$, we have
\begin{equation}
    T_3
    = \sum_{j=\lfloor j_* - \ell_n^{-} \rfloor}^{\lfloor j_* + \ell_n^{+} \rfloor}  \frac{1}{n^{1/4}}h_2\left(x\left(j/n\right)\right)
    = o(1),
\end{equation}
where the function $h_2$ is defined in~\eqref{eq:defh1}.
}
\begin{proof}
The proof will follow the same lines as the proof \autoref{ref:sumT2}. Let us start with the asymptotic of $h_2(x)$ for large $x$. Doing the change of variable $y\mapsto y/x$ in the integral expression of $f_1$ and $f_2$ \eqref{eq:f1} then expanding in powers of $x$ one finds that, as $x\to \infty$,
\Eq{}{
f_1(x)
=&\,
\frac{2\sqrt{\pi}}{\sqrt{\gamma_Q}}\frac{1}{x}
\left(
1+\frac{7}{6\gamma_Q}x^{-4}
+\frac{385}{72\gamma_Q^2}x^{-8}
+O(x^{-12})
\right)\\
f_2(x)=&\,-\frac{\sqrt{\pi}}{3\,\gamma_Q^{3/2}}\left(\frac{\gamma_Q}{r_*} + \upsilon_Q\right)
+\frac{\sqrt{\pi}}{3\,\gamma_Q^{5/2}}\left(3\,\upsilon_Q-\frac{16\,\gamma_Q}{r_*} \right)x^{-4}\\
&+\frac{7\sqrt{\pi}}{216\,\gamma_Q^{7/2}}\left(2201\,\upsilon_Q-\frac{1735\,\gamma_Q}{r_*}\right)x^{-8}
+ O(x^{-12})
}

Using those asymptotics we get that
\Eq{}{
h_2(x)=\frac{f_2(x)}{f_1(x)}
+
\frac{ \left(\upsilon_Q+\frac{\gamma_Q}{r_*}\right)}{6\gamma_Q}x
=
-\left(\frac{\upsilon_Q}{6\gamma_Q^2}+\frac{10}{3 r_*\gamma_Q}\right)x^{-3}
+O(x^{-7}),\quad x\to\infty.
}
Following the same steps as in the proof of \autoref{ref:sumT2} one arrives at
\Eq{}{
T_3= \sum_{j=\lfloor j_* - \ell_n^{-} \rfloor}^{\lfloor j_* + \ell_n^{+} \rfloor}  \frac{1}{n^{1/4}}h_2\left(x\left(j/n\right)\right)=\frac{r_*\gamma_Q}{4}\int_{-\xi_n}^{\xi_n}x^2h_2(x)\,dx+o(1).
}
Noticing that $x\mapsto x^2h_2(x)$ is an odd function, the above integral vanishes and this finishes the proof.
\end{proof}
\lem{ref:sumT4}{
As $n \to \infty$, we have
\Eq{}{
T_4:=\sum_{j \notin [\,\lfloor j_* - \ell_n^{-} \rfloor,\, \lfloor j_* + \ell_n^{+} \rfloor\,]} \frac{1}{n}\mathcal{B}(r_\tau)=&\,\mathcal{F}_{Q}[\mathbb{A}_{r_0,r_1}]+\frac{1}{3}\log(\frac{r_1}{r_0})-\frac{1}{4}\log\!\left(\frac{\Delta Q(r_1)}{\Delta Q(r_0)}\right)\\
&+\frac{7 r_*}{12}\frac{ n^{1/4}}{\xi_n}+o(1),
}
where the function $\mathcal{B}$ is defined in~\eqref{eq: def B} and $\mathcal{F}_{Q}[\mathbb{A}_{r_0,r_1}]$ is given in \eqref{eq: def FQ}.
}
\begin{proof}
Applying once again Euler-Maclaurin formula, the boundary terms are negligible and one gets
\Eq{}{
\sum_{j \notin [\,\lfloor j_* - \ell_n^{-} \rfloor,\, \lfloor j_* + \ell_n^{+} \rfloor\,]} \frac{1}{n}\mathcal{B}(r_\tau)=\left(\int_{r_0}^{r_*-\frac{\xi_n}{n^{1/4}}}+\int^{r_1}_{r_*+\frac{\xi_n}{n^{1/4}}}\right)\frac{1}{2}\mathcal{B}(r)r \Delta Q(r)dr+o(1).
}
One now wishes to extend the domain of integration of the above integral to the full interval $[r_0,r_1]$ while controlling the error. For that, we recall that, by assumption, $\Delta Q(r_*)=\partial_r\Delta Q(r_*)=0$ which implies that the function  $r\mapsto\mathcal{B}(r)r \Delta Q(r)$ has a double pole at $r=r_*$. More precisely, we have
\Eq{}{
\frac{1}{2}\mathcal{B}(r)r \Delta Q(r)=\frac{C_{1,Q}}{(r-r_*)^2}+\frac{C_{2,Q}}{(r-r_*)}+O(1),\qquad r\to r_*,
}
with $C_{1,Q}$ and $C_{2,Q}$ given in \eqref{eq: def CQ}. Hence, we obtain
\Eq{eq: T4 interm}{
T_4=& \int_{r_0}^{r_1}\left[ \frac{1}{2}\mathcal{B}(r)r \Delta Q(r)-\frac{C_{1,Q}}{(r-r_*)^2}-\frac{C_{2,Q}}{(r-r_*)}\right]dr\\
&+\left(\int_{r_0}^{r_*-\frac{\xi_n}{n^{1/4}}}+\int^{r_1}_{r_*+\frac{\xi_n}{n^{1/4}}}\right)\left[\frac{C_{1,Q}}{(r-r_*)^2}+\frac{C_{2,Q}}{(r-r_*)}\right]dr +o(1)
}
where the last integral can be computed explicitly
\Eq{eq: int correction}{
\left(\int_{r_0}^{r_*-\frac{\xi_n}{n^{1/4}}}+\int^{r_1}_{r_*+\frac{\xi_n}{n^{1/4}}}\right)\left[\frac{C_{1,Q}}{(r-r_*)^2}+\frac{C_{2,Q}}{(r-r_*)}\right]dr=&\,\frac{2C_{1,Q} n^{1/4}}{\xi_n}
-\frac{C_{1,Q}}{r_1-r_*}
-\frac{C_{1,Q}}{r_*-r_0}
\\
&+C_{2,Q} \log\!\left(\abs{\frac{r_1-r_*}{r_*-r_0}}\right).
}
Using the explicit expression of $\mathcal{B}$ \eqref{eq: def B}, the first integral in \eqref{eq: T4 interm} can be put in the form
\Eq{eq: IBP}{
&\int_{r_0}^{r_1}\left[- \frac{1}{16} r\frac{\partial_r^2 \Delta Q(r)}{\Delta Q(r)}
- \frac{19}{48 } \frac{\partial_r \Delta Q(r)}{\Delta Q(r)}
+ \frac{5}{48}r \left(\frac{\partial_r \Delta Q(r)}{\Delta Q(r)}\right)^2
+ \frac{1}{6r} -\frac{C_{1,Q}}{(r-r_*)^2}-\frac{C_{2,Q}}{(r-r_*)}\right]dr\\
=&\int_{r_0}^{r_1}\left[-\frac{7}{24}r \left(\frac{\partial_r \Delta Q(r)}{\Delta Q(r)}\right)^2
+\frac{1}{3 } r\frac{\partial_r^2 \Delta Q(r)}{\Delta Q(r)}
-\frac{r C_{2,Q}+C_{1,Q}}{(r-r_*)^2}\right]dr+\frac{1}{6}\log(\frac{r_1}{r_0})\\
&- \frac{19}{48 }\left[ r_1\frac{\partial_r \Delta Q(r_1)}{\Delta Q(r_1)}-r_0\frac{\partial_r \Delta Q(r_0)}{\Delta Q(r_0)}\right]-C_{2,Q}\left[\frac{r_1}{r_1-r_*}
+\frac{r_0}{r_*-r_0}\right],
}
where we have used the following integration by parts
\Eq{}{
& \int_{r_0}^{r_1}\left[- \frac{19}{48 } \frac{\partial_r \Delta Q(r)}{\Delta Q(r)}
-\frac{C_{2,Q}}{(r-r_*)}\right]dr\\
=&\,- \frac{19}{48 }\left[ r_1\frac{\partial_r \Delta Q(r_1)}{\Delta Q(r_1)}-r_0\frac{\partial_r \Delta Q(r_0)}{\Delta Q(r_0)}\right]-C_{2,Q}\left[\frac{r_1}{r_1-r_*}
+\frac{r_0}{r_*-r_0}\right]\\
&+\int_{r_0}^{r_1}r\left[\frac{19}{48 }\left( \frac{\partial_r^2 \Delta Q(r)}{\Delta Q(r)}-\left(\frac{\partial_r \Delta Q(r)}{\Delta Q(r)}\right)^2\right)
-\frac{C_{2,Q}}{(r-r_*)^2}\right]dr.
}
Combining \eqref{eq: IBP} and \eqref{eq: int correction} finishes the proof.
\end{proof}

Combining \autoref{ref:sumT1}, \autoref{ref:sumT2}, \autoref{ref:sumT3} and \autoref{ref:sumT4} gives the claim of  \autoref{theo1} in the annulus case, i.e. $r_0>0$. Let us now treat the case $r_0=0$. 

\subsection{Disk case}
Let us denote $Z_{n,\mathrm{disk}}$ to be the partition function in the case the droplet is a disk, i.e. $r_0=0$. In this setting, the asymptotic behaviour of $u_j$ differs from the case $r_0>0$ for small indices and it is given by \autoref{lem: uj asym6}. 
In \cite{AFLS25}, it was remarked that this same change in asymptotic behaviours occurs in the study \cite{BKS23}
of the partition function $Z_n^{\mathrm{sh}}$ (which stands for the partition function with strictly subharmonic potential).
Since this common alteration is independent of the subharmonicity property of the potential away from the origin, we have 
\begin{equation}\label{Z1}
  \log(Z_{n,\mathrm{disk}}^{\rm sh})
  - \log(Z_{n,\mathrm{annulus}}^{\rm sh})
  \;=\;
  \log(Z_{n,\mathrm{disk}})
  - \log(Z_{n,\mathrm{annulus}}) .
\end{equation}
Using \eqref{eq: functional reduce}, the left-hand side of \eqref{Z1} can be read off from the results of \cite[Theorem 1.1 (i) and Theorem 1.2 (i)]{BKS23}, and is equal to
\Eq{}{
  -\frac{1}{12}\,\log(n) + \zeta'(-1) + \mathcal{F}_{Q}[\mathbb{D}_{r_1}] - \mathcal{F}_{Q}[\mathbb{A}_{r_0,r_1}],
}
thus completing the proof of \autoref{theo1}.
\end{proof}

\section{Local correlations}\label{local correlations}
\subsection{Local Kernel estimate}\label{localkernelestimateorg}

The main goal of this section is to estimate the range of indices in the sum \eqref{eq: def kernel} having a non-negligible contribution. The modulus of the kernel is  bounded by
\Eq{}{
\abs{K\left(z_1,z_2 \right)}\leq e^{-\frac{n}{2}(Q(z_1)+Q(z_2))}\sum_{j=0}^{n-1}\frac{\abs{p_{j,n}(z_1)} \abs{p_{j,n}(z_2)}}{\norm{p_{j,n}}^2}=e^{-\frac{n}{2}(Q(z_1)+Q(z_2))}\sum_{j=0}^{n-1}\frac{\abs{z_1}^j \abs{z_2}^j}{\norm{p_{j,n}}^2}
}
Going into polar coordinates $z_1=\varrho_1e^{i\theta_1}$ and $z_2=\varrho_2e^{i\theta_2}$ and considering the radially symmetric potential $Q$ \eqref{eq:complex potential}, it becomes
\Eq{}{
\abs{K\left(\varrho_1e^{i\theta_1},\varrho_2e^{i\theta_2} \right)} \leq \sum_{j=0}^{n-1}\frac{\varrho_1^{j} e^{-\frac{n}{2}Q(\varrho_1)} }{\norm{p_{j,n}}} \frac{\varrho_2^{j} e^{-\frac{n}{2}Q(\varrho_2)} }{\norm{p_{j,n}}}.
}
If we are interested in values of radii $\varrho_1$ and $\varrho_2$ which are around $r_*$, the radii for which the gap starts to appear, then the following proposition gives a bound on the summands.
\prop{kernelestimate}{
    Let $z\in\mathbb{C}$ be such that $||z|-r_*|\leq\frac{\log(n)}{n^{1/4}}$. For $j\notin  [j_*-\ell_n^{-},j_*+\ell_n^{+}]$, 
    \begin{equation}\label{eq:keres}
        \frac{|p_{j,n}(z)|^2}{{\norm{p_{j,n}}^2}}\lesssim n e^{-c \xi_n^4 },
    \end{equation}
    with $\xi_n$ given in \eqref{eq: def xi} and $c$ a positive constant independent of $n$.
}
\begin{proof}
    A standard pointwise $L^2$ estimate  from potential theory \cite[Lemma 2.4]{A21} gives
    \begin{equation}
         \frac{1}{n}e^{-n V_\tau(r_\tau)}\lesssim {\norm{p_{j,n}}^2}.
    \end{equation}
with $\lesssim$ as in \autoref{nota: leq geq}. Thus, going into polar coordinates $z=re^{i\theta}$, we have 
\begin{equation}\label{eq:sumbound}
    \frac{|p_{j,n}(z)|^2}{{\norm{p_{j,n}}^2}}\lesssim n e^{-n (V_\tau(r)-V_\tau(r_\tau))}.
\end{equation}
The goal is then to estimate the difference
\begin{equation}\label{eq: diff V}
    V_\tau(r)-V_\tau(r_\tau)=V_\tau(r)-V_\tau(r_*)+ V_\tau(r_*)- V_\tau(r_\tau)\geq 0.
\end{equation}
Let us note that this difference is always positive as $r_\tau$ is the unique minimum of $V_\tau$. While the term $V_\tau(r)-V_\tau(r_*)$ is guaranteed to be small under the assumption that $|r-r_*|\leq \frac{\log(n)}{n^{1/4}}$, the aim is to determine how small can $V_\tau(r_*)- V_\tau(r_\tau)$ be under the assumption that $j\notin  [j_*-\ell_n^{-},j_*+\ell_n^{+}]$. 

Let us introduce for convenience the two functions
\begin{equation}
  A(r):=V_\tau(r)-V_\tau(r_*),\qquad B(\tau):= V_\tau(r_*)-V_\tau(r_\tau).  
\end{equation}
Focusing first on $B$, we have
\begin{equation}
   B'(\tau)=2\log (r_\tau)-2\log (r_*).
\end{equation}
Reminding that $\tau\in[0,\frac{n-1}{n}]\setminus [\tau_*-\frac{\ell_n^{-}}{n},\tau_*+\frac{\ell_n^{+}}{n}]$ and $r_\tau\in [r_0,r_1]$, due to the compactness of their range, it is always possible to find some $n$ independent constants such that 
\begin{equation}
   \abs{B'(\tau)}=2  \abs{\log(1+\frac{r_\tau-r_*}{r_*})}\gtrsim |r_\tau-r_*|\gtrsim |\tau-\tau_*|^{1/3}
\end{equation}
where we have made use of the expansion \eqref{eq:locexp}. Thus, also by compactness of the range of $\tau$, we have
\Eq{eq: B bound}{
|B(\tau)|\gtrsim\abs{\tau-\tau_*}^{4/3}.
}

Moving to $A(r)$, a simple Taylor expansion around $r=r_*$ gives the upper bound
\Eq{eq: A bound}{
\abs{A(r)}= \abs{V_\tau(r)-V_\tau(r_*)}\leq |V'_\tau(r_*)||r-r_*|&+\frac{|V^{(2)}_\tau(r_*)|}{2}|r-r_*|^2+\frac{|V^{(3)}_\tau(r_*)|}{3!}|r-r_*|^3\\&+\sup_{r_0<|\varrho|<r_1} \frac{|V^{(3)}_\tau(\varrho)|}{4!}|r-r_*|^4.
}
One can then express $V'_\tau(r_*)$, $V^{(2)}_\tau(r_*)$ and $V^{(3)}_\tau(r_*)$ in terms of $\tau-\tau_*$. For the first derivative we use \eqref{eq: tau*} to arrive at
\Eq{}{
V'_\tau(r_*)=q'(r_*) - \frac{2\tau}{r_*}=-2\,\frac{\tau-\tau_*}{r_*}.
}
Then, following from \eqref{eq: deriv V} we get
\Eq{}{
|V'_\tau(r_*)|=2\,\frac{\abs{\tau-\tau_*}}{r_*},\qquad |V^{(2)}_\tau(r_*)|=2\,\frac{\abs{\tau-\tau_*}}{r_*^2}, \qquad|V^{(3)}_\tau(r_*)|=4\,\frac{\abs{\tau-\tau_*}}{r_*^3}.
}
We are now almost ready to find a lower bound for \eqref{eq: diff V}. First we denote $C:=\sup_{r_0<|\varrho|<r_1} |V^{(4)}_\tau (\varrho)|/(4!)$ and define 
\Eq{eq: def elln}{
\ell_n:=\min\{\ell_n^-,\ell_n^+ \},
}
which is an increasing sequence behaving as $\ell_n\sim \frac{\varkappa_Q}{3!} n^{1/4}\xi_n^3$, as $n\to\infty$, due to \eqref{eq:locexp} and \eqref{eq:defdn}, where $\xi_n$ is given in \eqref{eq: def xi} and $\varkappa_Q$ is given in \eqref{eq:vcon}.
Then, by reverse triangle inequality and using the assumptions $|r-r_*|<\frac{\log(n)}{n^{1/4}}$, $\tau\notin [\tau_*-\frac{\ell_n^{-}}{n},\tau_*+\frac{\ell_n^{+}}{n}]$, the bounds \eqref{eq: B bound} and \eqref{eq: A bound} leads to
\Eq{eq: ineq}{
    n(V_\tau(r)-V_\tau(r_\tau))\geq n\abs{|B(\tau)|-|A(r)|}\geq F_n,
}
where $F_n$ is the positive sequence
\Eq{eq: F_n}{
F_n:=&\,\Big|\alpha\, n|\tau-\tau_*|^{4/3}-C \log(n)^4\\&-\left(\frac{2}{r_*}n^{3/4}\log(n)+\frac{1}{r_*^2}n^{1/2}\log(n)^2+\,\frac{2}{3 r_*^3}n^{1/4}\log(n)^3\right)|\tau-\tau_*|\Big|,
}
with $\alpha>0$, independent of $n$. 

We now wish to show that the sequence $F_n$ cannot be too small, as $n\to\infty$, when imposing  $\tau\notin [\tau_*-\frac{\ell_n^{-}}{n},\tau_*+\frac{\ell_n^{+}}{n}]$. In particular we want to show that it is bounded from below by $c \log(n)^4$ for some constant $c > 0$ independent of $n$. As $\abs{\tau-\tau_*}$ is implicitly a function of $n$ we introduce the notation
\Eq{}{
L_n:=n\abs{\tau-\tau_*}
}
and proceed with a case distinction depending on the large $n$ asymptotic behavior of the sequence $L_n$. From the definition of $\ell_n$ \eqref{eq: def elln} and of \eqref{eq:defdn} we have that $L_n\geq\ell_n>0$. We reminded the asymptotic behavior of $\ell_n$ below \eqref{eq: def elln}. If $L_n$ is of the same order as $\ell_n$, some cancellation of leading term can happen in the expression of $F_n$ \eqref{eq: F_n}, hence the following distinction:
\begin{itemize}
    \item If $L_n/\ell_n\underset{n\to\infty}{\to}\ell<\infty$, then $F_n$ has the following asymptotic behavior
\Eq{}{
F_n\sim \alpha\left(\frac{\ell \varkappa_Q}{3!}\right)^{4/3}\xi_n^4, \quad n\to\infty.
}
\item If $L_n/\ell_n\underset{n\to\infty}{\to}\infty$, then 
\Eq{}{
F_n\sim \alpha\left(\frac{ \varkappa_Q}{3!}\right)^{4/3}\xi_n^4 \left(\frac{L_n}{\ell_n}\right)^{4/3}, \quad n\to\infty,
}
which is clearly larger than $\xi_n^4$ for $n$ large enough.
\end{itemize}
This finishes to show that, for sufficiently large $n$, we have
\begin{equation}\label{eq:es1}
    n(V_\tau(r)-V_\tau(r_\tau))\geq c\xi_n^4.
\end{equation}
for some $c>0$, independent of $n$. Plugging \eqref{eq:es1} in \eqref{eq:sumbound} finishes the proof.
\end{proof}

\autoref{kernelestimate} motivates the definition of a local kernel, near the radius of the gap formation,
   \begin{equation}\label{eq: Kn gap}
       K_n^{\rm g}(z_1,z_2):=\sum_{j=\lfloor j_* - \ell_n^{-} \rfloor}^{\lfloor j_* +\ell_n^{+} \rfloor}\frac{p_{j,n}(z_1)\overline{p_{j,n}(z_2)}}{{\norm{p_{j,n}}^2}}.
   \end{equation}
As a consequence of \autoref{kernelestimate} we obtain the following corollary.
\coro{coro:tail}{
 Let $z_1,z_2\in \mathbb{C}$. If $||z_1|-r_*|,||z_2|-r_*|\leq \frac{\log (n)}{n^{1/4}}$  then, as $n\to\infty$, the kernel admits the following expansion
    \begin{equation}
        K_{n}(z_1,z_2)=K_n^{\rm g}(z_1,z_2)+O(n^2 e^{-c\xi_n^4}),
    \end{equation}
    for some $c>0$ independent of $n$ and $\xi_n$ given in \eqref{eq: def xi}.
}

\subsection{Proof of \autoref{theo2} }\label{sec:theo2}
Before starting the proof of \autoref{theo2}, let us state the following corollary, which follows directly from \autoref{prop:norm}.
\coro{cor:norm}{
For $j\in  [j_*-\ell_n^{-},j_*+\ell_n^{+}]$, we have, as $n\rightarrow\infty$,
\begin{equation}
    \frac{1}{{\norm{p_{j,n}}^2}}=n^{1/4}e^{n V_\tau(r_\tau)}\frac{1}{2r_* f_1(x)}\left(1-\frac{1}{n^{1/4}}\frac{f_2(x)}{f_1(x)}+O(n^{-1/2})\right),
\end{equation}
where $f_1$ and $f_2$ are given in \eqref{eq:f1}, and $x = n^{1/4} (r_\tau-r_*)$.
}

\begin{proof}[Proof of \autoref{theo2}]
Let $z_1,z_2$ be such as in \autoref{theo2} and denote 
\Eq{eq: varrho}{
 \varrho_k:=\abs{z_k}=\abs{r_*+\frac{\xi_k}{(\gamma_Q\, n)^{1/4}}},\quad \theta_k:=\arg(z_k),\quad k=1,2.
}
We will make use of the following expansions, as $n\to \infty$,
\Eq{eq: expan rho_k}{
\varrho_k=r_*+\frac{\Re(\xi_k)}{(\gamma_Q\, n)^{1/4}}+O(n^{-1/2}),\qquad
\theta_k= \frac{\Im(\xi_k)}{r_*(\gamma_Q\, n)^{1/4}}+O(n^{-1/2}).
}
Rewriting the polynomials
\begin{equation}\label{eq:pj}
    p_{j,n}(z_k)=e^{ij\theta_k}e^{-\frac{n}{2} V_\tau(\varrho_k)},\quad k=1,2,
\end{equation}
and using \autoref{cor:norm} it follows that
\Eq{eq: summand pj}{
    \frac{p_{j,n}(z_1)\overline{p_{j,n}(z_2)}}{{\norm{p_{j,n}}^2}}=&\,n^{1/4}\exp(i \,j (\theta_1-\theta_2))e^{-\frac{n}{2}(V_\tau(\varrho_1)-V_\tau(r_\tau))}e^{-\frac{n}{2}(V_\tau(\varrho_2)-V_\tau(r_\tau))}\\
    &\times\frac{1}{2r_* f_1(x)}\left(1-\frac{1}{n^{1/4}}\frac{f_2(x)}{f_1(x)}+O(n^{-1/2})\right).
}

One can then Taylor expand the function $r\mapsto V_\tau(r)-V_\tau(r_\tau)$ around $r=r_\tau$. Evaluating the expansion for $r=\varrho_k$ yields
\Eq{}{
 V_\tau(\varrho_k)-V_\tau(r_\tau)=&\,\frac{1}{2}V_\tau''(r_\tau)(\varrho_k-r_\tau)^2+\frac{1}{3!}V_\tau^{(3)}(r_\tau)(\varrho_k-r_\tau)^3\\
 &+\frac{1}{4!}V_\tau^{(4)}(r_\tau)(\varrho_k-r_\tau)^4+\frac{1}{5!}V_\tau^{(5)}(r_\tau)(\varrho_k-r_\tau)^5+O((\varrho_k-r_\tau)^6).
}
The above expansion is a function of $r_\tau$ which can be further expanded around $r_\tau=r_*$ using \eqref{eq:vder} and 
\Eq{}{
\delta_r =r_\tau-r_*= \frac{v}{(\gamma_Q n)^{1/4}}.
}
Combining this with the expansion \eqref{eq: expan rho_k}, one obtains a similar expansion to \eqref{eq: phase expansion} where $y$ is replaced by $(\Re(\xi_k)-v)/\gamma_Q^{1/4}$. It then leads to the following expansion, as $n \to \infty$, 
\Eq{eq: exp exp}{
&\exp\left[-n( V_\tau(\varrho_k)-V_\tau(r_\tau))\right]\\
=&\exp[-v^4 \mathbf{P}'\left(\frac{\Re(\xi_k)}{v}-1\right)]\Bigg(1-\frac{1}{n^{1/4}}\frac{\Phi_Q(\xi_k,v)}{\gamma_Q^{1/4}}+O(n^{-1/2} )\Bigg)
}
where $\mathbf{P}$ is given by \eqref{eq: def bf P} and 
\Eq{eq: def phiQ}{
\Phi_Q(\xi_k,v):=\left(\frac{r_*\upsilon_Q}{\gamma_Q}+1\right)\frac{v^5}{2(3!)}\left(\frac{\Re(\xi_k)}{v}-1\right)^2+\frac{r_*\upsilon_Q}{\gamma_Q}v^5 \mathbf{P}\left(\frac{\Re(\xi_k)}{v}-1\right).
}
As $v=\gamma_Q^{1/4}x$, plugging back \eqref{eq: exp exp} into \eqref{eq: summand pj} yields
\Eq{eq:kersum}{
 \frac{p_{j,n}(z_1)\overline{p_{j,n}(z_2)}}{{\norm{p_{j,n}}^2}}=&\,\frac{n^{1/4}}{2r_*}\exp(i \,j (\theta_1-\theta_2))\frac{\exp\left(-\frac{v^4}{2}\left[ \mathbf{P}'\left(\frac{\Re(\xi_1)}{v}-1\right)+ \mathbf{P}'\left(\frac{\Re(\xi_2)}{v}-1\right)\right]\right)}{f_1\left(\gamma_Q^{-1/4} v\right)}\\
 &\times\left(1-\frac{1}{n^{1/4}}\left[\frac{f_2\left(\gamma_Q^{-1/4} v\right)}{f_1\left(\gamma_Q^{-1/4} v\right)}+\frac{\Phi_Q(\xi_1,v)+\Phi_Q(\xi_2,v)}{2r_*\gamma_Q^{1/4}}\right]+O(n^{-1/2})\right),
}
where, after performing the change of variable $y\mapsto (v/\gamma_Q^{1/4})y$ in the integral expression of $f_1$ and $f_2$  \eqref{eq:f1}, one has the explicit expression
\Eq{}{
\frac{f_2\left(\gamma_Q^{-1/4} v\right)}{f_1\left(\gamma_Q^{-1/4} v\right)}=\frac{1}{r_*\gamma_Q^{1/4}\mathcal{P}(v^4)}\int_{-\infty}^{\infty}\left[v y-\Phi_Q(vy+1,v)\right]\exp[-v^4\mathbf{P}'\left(y\right)]dy.
}
Using the change of variable \eqref{eq:defp}, where we recall that $r_\tau - r_*=v/(\gamma_Q\,n)^{1/4}$, one has 
\Eq{}{
j=j_*+n\, p\left(\frac{v}{(\gamma_Q\,n)^{1/4}}\right)
}
and, as $n\to\infty$,
\Eq{eq: p expansion}{
  p\left(\frac{v}{(\gamma_Q\,n)^{1/4}}\right)&=\frac{r_* \gamma_Q^{1/4}}{2\,n^{3/4}}\frac{v^3}{3!}+\left(\frac{r_*\upsilon_Q}{\gamma_Q} + 4 \right)\frac{1}{2\, n} \frac{v^4}{4!}+ O(n^{-5/4}),\\
  p'\!\left(\frac{v}{(\gamma_Q\,n)^{1/4}}\right)&=\frac{r_*\gamma_Q^{1/2} }{2 \,n^{1/2}}\frac{v^2}{2} + \left(\frac{r_*\upsilon_Q}{\gamma_Q^{3/4}} + 4 \gamma_Q^{1/4}\right)\frac{1}{2\, n^{3/4}} \frac{v^3}{3!}
    + O(n^{-1})\\
    &=\frac{r_*\gamma_Q^{1/2} }{2 \,n^{1/2}}\frac{v^2}{2}\left(1+\frac{1}{r_*(\gamma_Q n)^{1/4}}\left(\frac{r_*\upsilon_Q}{\gamma_Q} + 4\right)\frac{v}{3} +O(n^{-1/2})\right)
  }
Combining the above expansion with the expansion \eqref{eq: expan rho_k}, we find that, as $n\to \infty$, the complex phase $j (\theta_1-\theta_2)$ in \eqref{eq:kersum} satisfies
  \Eq{eq: exp phase}{
  j (\theta_1-\theta_2)=j_*\frac{\Im(\xi_1-\xi_2)}{r_*(\gamma_Q\, n)^{1/4}}+\Im(\xi_1-\xi_2)\frac{v^3}{12}+\frac{\Im(\xi_1-\xi_2)}{r_*(\gamma_Q\, n)^{1/4}}\left(\frac{r_*\upsilon_Q}{\gamma_Q} + 4\right)\frac{1}{2} \frac{v^4}{4!}+O(n^{-1/2}).
  }
With all the required asymptotic expansions in hand, we proceed as in the proof of \autoref{ref:sumT2}, applying the Euler--Maclaurin formula (cf. \autoref{prop: EM formula}) to \eqref{eq:kersum} together with \eqref{eq: exp phase} and \eqref{eq: p expansion} to obtain
\Eq{eq: exp kernel 2}{
  &K_n^{\rm g}(z_1,z_2)=\sum_{j=\lfloor j_* - \ell_n^{-} \rfloor}^{\lfloor j_* + \ell_n^{+} \rfloor}  \frac{p_{j,n}(z_1)\overline{p_{j,n}(z_2)}}{{\norm{p_{j,n}}^2}}\\
& \qquad =\frac{\sqrt{\gamma_Q\, n} }{8}e^{i \,j_*\frac{\Im(\xi_1-\xi_2)}{r_*(\gamma_Q\, n)^{1/4}}}\int_{-\infty}^{\infty}dv\, e^{i\Im(\xi_1-\xi_2)\frac{v^3}{12}} \frac{\exp\left(-\frac{v^4}{2} \sum_{k=1,2}\mathbf{P}'\left(\frac{\Re(\xi_k)}{v}-1\right)\right)}{\mathcal{P}(v^4)}\abs{v}\\
& \qquad \quad\times\Bigg(1-\frac{1}{r_*(\gamma_Q n)^{1/4}}\Bigg[\int_{-\infty}^{\infty}\left[v y-\Phi_Q(vy+1,v)\right]\frac{\exp[-v^4\mathbf{P}'\left(y\right)]}{\mathcal{P}(v^4)}dy+\sum_{k=1,2}\frac{\Phi_Q(\xi_k,v)}{2}\\
& \qquad \quad-\left(\frac{r_*\upsilon_Q}{\gamma_Q} + 4\right)\frac{v}{3}-i\Im(\xi_1-\xi_2)\left(\frac{r_*\upsilon_Q}{\gamma_Q} + 4\right)\frac{1}{2} \frac{v^4}{4!}\Bigg]+O(n^{-1/2})\Bigg).
}
Let us point out that the factor 
\Eq{}{
e^{i \,j_*\frac{\Im(\xi_1-\xi_2)}{r_*(\gamma_Q\, n)^{1/4}}}=e^{i \,j_*\frac{\Im(\xi_1)}{r_*(\gamma_Q\, n)^{1/4}}}/e^{i \,j_*\frac{\Im(\xi_2)}{r_*(\gamma_Q\, n)^{1/4}}}
}
is a cocycle and can be eliminated by acting on the rows and columns of the determinant \eqref{eq: dpp}. Therefore, as the kernel is defined up to a cocycle, one can ignore it. Appealing to \autoref{coro:tail} then finishes the proof.
\end{proof}

\subsection{Mean level spacing}
\begin{proof}[Proof of \autoref{cor:mls}]
We remind that by \autoref{def: mls} the mean level spacing $s_n$ is such that 
\Eq{}{
\mathbb{E}\left[N_n\!\bigl(\mathbb{D}(r_*,s_n)\bigr)\right]=1.
}
Let us start from the left-hand side of the above expression. Since the macroscopic probability density is given by $z\mapsto K_n(z,z)/n$, the average number of particle in a disk of radius $s_n$ centered around $r_*$ is then $n$ times the probability to find at least one particle in this disk, i.e.
\begin{equation}\label{eq:ENn}
     \mathbb{E}\left[N_n\!\bigl(\mathbb{D}(r_*,s_n)\bigr)\right]
    =\int_{|z-r_*|<s_n} K_n(z,z)\,dA(z).
\end{equation}
For $z$ in a small neighbourhood of $r_*$ the first order approximation of the kernel is given by $K_n(z,z)=K_n(r_*,r_*)+o(1)$, as $z\to r_*$. The leading term asymptotic of the $n$-dependent constant $K_n(r_*,r_*)$ is then given by \autoref{coro2} for $\xi=0$, i.e.
\Eq{}{
K_n(r_*,r_*)=\frac{\sqrt{\gamma_Q\, n} }{4}\rho(0)+o(1),\qquad n\to\infty,
}
with $\rho$ defined in \eqref{eq: def rho}. Multiplying by the area of the disk $\pi s_n^2$ yields
\begin{equation}
     \mathbb{E}\left[N_n\!\bigl(\mathbb{D}(r_*,s_n)\bigr)\right]=1
    =\pi\frac{\sqrt{\gamma_Q\, n} }{4}\rho(0)\,s_n^2+o(1),\qquad n\to\infty.
\end{equation}
which then implies
\begin{equation}
 s_n=\frac{2}{\sqrt{\pi\rho(0)}}\,\frac{1}{(\gamma_Q\, n)^{1/4}}+ o\left(\frac{1}{n^{1/4}} \right),\qquad n\to\infty,
\end{equation}
and finishes the proof.
\end{proof}
\section{Discussion and outlook}\label{discussion}
In this work, we studied a 2D Coulomb gas at inverse temperature $\beta=2$ and considered a family of radially symmetric potentials for which a gap is about to appear inside the droplet---which is either a disk or an annulus---around a circle of radius $r_*>0$. 
In particular, we studied the impact of this critical phenomena, both on the macroscopic level via the large $n$ asymptotic of the partition function, cf. \autoref{theo1}, and on the microscopic level via the local correlation functions which are controlled, in this case, by a new universal kernel, cf. \autoref{theo2}.

Several questions remain open. First, in this work, we restricted ourselves to radially symmetric potentials. If the radial symmetry is broken, it remains to be understood whether the universality of the new kernel in \eqref{eq: double scaling limit} persists. One can conjecture that it does, as long as the gap appears along a smooth enough curve, since the limiting kernel captures only local correlations. Regarding the expansion of the partition function, without radial symmetry, the expansion is expected to still have a term of order $n^{1/4}$, reminiscent of the criticality of the phenomena, but the associated constant is expected to be different.
 
It would also be of great interest to study the fluctuations of linear statistics of the point process \eqref{eq:2D CG}, near the critical radius $r_*$. See \cite{F23} for a survey of results in this direction. 

We also note that, as a generalisation of \autoref{example potential}, one can consider a one parameter family of potentials defined by the difference  
\begin{equation}\label{eq:cpot}
    Q_t(z)=Q_1(z)-tQ_2(z),
\end{equation}
with $Q_1$ and $Q_2$ strictly subharmonic and such that  
\begin{equation}
    \Delta Q_1(z)\geq \Delta Q_2(z),
\end{equation}
with equality for $|z|=r_*$, while strict inequality for $|z|\neq r_*$.  As $t$ varies in some small interval centered around $t=1$, the number of connected components of the support of the equilibrium measure increases. The present work concerns the case $t=1$. A natural next step is to study the point process in the scaling limit where $\abs{t-1}=c_n$ for some sequence $c_n\to 0$. We expect in this regime an interpolation between the Ginibre bulk kernel and the critical kernel in \eqref{eq: def kernel}. This will be the object of another article. 

Finally, our large $n$ expansion of the partition function reveals a new term of order $n^{1/4}$. This scale is inversely proportional to the mean level spacing of particles in the vicinity of the critical phenomena---near emergence of a gap---and dependent on the degree of criticality of the potential. We give \autoref{conj: Z_n} and \autoref{conj: kernel m} regarding the expansion of the partition function and the microscopic correlation functions under higher degree of criticality of the potential. 
We plan to return to these questions in future works.

\appendix
\section{Euler-Maclaurin formula}\label{appendix}
The Euler-Maclaurin summation formula \cite{DLMF}, which we state below for reference, plays a crucial role in our analysis by allowing us to derive the large-$n$ asymptotic behavior of the sums considered.
\prop{prop: EM formula}{(Euler-Maclaurin formula) Let $f(x)$ be $2k$ times differentiable  on the interval $[p,q]$. Denote by $\{B_{2k} \}$ the even indexed Bernoulli numbers ($B_2=\frac{1 }{ 6}, B_4 = - \frac{1 }{ 30}, \dots$)
\Eq{}{
\sum _{i=p+1}^{q-1}f(i)=\int _{p}^{q}f(x)~{\rm {d}}x-{\frac {f(p)+f(q)}{2}}+\sum _{j=1}^{k}{\frac {B_{2j}}{(2j)!}}\left(f^{(2j-1)}(q)-f^{(2j-1)}(p)\right)+R_{2k}
}
\Eq{}{
\sum _{i=p}^{q}f(i)=\int _{p}^{q}f(x)~{\rm {d}}x+{\frac {f(p)+f(q)}{2}}+\sum _{j=1}^{k}{\frac {B_{2j}}{(2j)!}}\left(f^{(2j-1)}(q)-f^{(2j-1)}(p)\right)+R_{2k},
}
where the remainder $R_{2k}$ satisfies the bound
$|R_{2k}| \le c_{2k}
\int_p^q | f^{(2k)}(x) | \, dx$, with $c_{2k} = 2 \zeta(2k) / (2 \pi)^{2k}$.
}

\section{First correction term for the kernel}\label{app correction term}
The proof of \autoref{theo2} contains the explicit expression of the first correction to the limiting kernel. However, since the correction term is cumbersome and offers little additional insight, we include it here for completeness. To make the structure apparent and for the sake of compactness, we introduce the following polynomial
\Eq{eq: def L}{
L(y):=\frac{y^2}{2(3!)}+\frac{y^3}{2(3!)}+\frac{y^4}{4!}+\frac{y^5}{5!},
}
as well as the functional $\mathscr{P}$ which acts on a polynomial function $\wp$ as
\Eq{eq: fct}{
\mathscr{P}[\wp](\theta):=\int_{-\infty}^{\infty}\wp(y)\exp[-\theta P\left(y\right)]dy,
}
where $P$ is given in \eqref{eq: pearcey}. We will denote
\Eq{}{
\mathcal{P}^{(k)}(\theta):=\mathscr{P}[x\mapsto x^k](\theta),\quad k\in\mathbb{N},
}
with $\mathcal{P}^{(0)}(\theta)\equiv \mathcal{P}(\theta)$. 
The first correction term of the limiting kernel can then be expressed in terms of the following functions, which are independent of the potential,
\Eq{eq: def Xi1}{
\Xi_*^{(1)}(\xi_1,\xi_2):=&\int_{-\infty}^{\infty}dv\, e^{i\Im(\xi_1-\xi_2)\frac{v^3}{12}} \frac{\exp\left(-\frac{v^4}{2}\left[ P\left(\frac{\Re(\xi_1)}{v}-1\right)+ P\left(\frac{\Re(\xi_2)}{v}-1\right)\right]\right)}{\mathcal{P}(v^4)}\abs{v}\\
  &\times\Bigg(-v\frac{\mathcal{P}^{(1)}(v^4)}{\mathcal{P}(v^4)}+\frac{v^5}{4!}\left[2\frac{\mathcal{P}^{(2)}(v^4)}{\mathcal{P}(v^4)} -\left(\frac{\Re(\xi_1)}{v}-1\right)^2-\left(\frac{\Re(\xi_2)}{v}-1\right)^2 \right]\\
  &+8\left[\frac{v}{3!}+i\Im(\xi_1-\xi_2)\frac{v^4}{4!}\right]\Bigg)
}
and
\Eq{eq: def Xi2}{
\Xi_*^{(2)}(\xi_1,\xi_2):=&\int_{-\infty}^{\infty}dv\, e^{i\Im(\xi_1-\xi_2)\frac{v^3}{12}} \frac{\exp\left(-\frac{v^4}{2}\left[ P\left(\frac{\Re(\xi_1)}{v}-1\right)+ P\left(\frac{\Re(\xi_2)}{v}-1\right)\right]\right)}{\mathcal{P}(v^4)}\abs{v}\\
  &\times\Bigg(\frac{v^5}{2}\left[2\frac{\mathscr{P}[L](v^4)}{\mathcal{P}(v^4)} -L\left(\frac{\Re(\xi_1)}{v}-1\right)-L\left(\frac{\Re(\xi_2)}{v}-1\right) \right]\\
  &+2\left[\frac{v}{3!}+i\Im(\xi_1-\xi_2)\frac{v^4}{4!}\right]\Bigg).
}
A stronger version of \autoref{theo2} is then given by the following theorem.
\theo{theo2 strong}{
For $k=1,2$, let $z_k=r_*+\frac{\xi_k}{(\gamma_Q\, n)^{1/4}}\in \mathbb C$, with $\gamma_Q$ given in \eqref{eq:vcon} and $|\xi_k|\leq \log(n)$. Under the same assumptions as \autoref{theo1}, the kernel \eqref{eq: def kernel}, up to a cocycle, admits the following uniform expansion, as $n \to \infty$,
\Eq{eq: kernel strong}{
 \frac{1}{\sqrt{\gamma_Q\,n}} K_{n}(z_1,z_2)=&\,\frac{1}{4}K_*(\xi_1,\xi_2)+\frac{1}{4 \varkappa_Q n^{1/4}}\left(\Xi_*^{(1)}(\xi_1,\xi_2)+\frac{r_*\upsilon_Q}{\gamma_Q} \Xi_*^{(2)}(\xi_1,\xi_2)\right)+O(n^{-1/2}),
}
where $K_*$, $\Xi_*^{(1)}$, $\Xi_*^{(1)}$, $\varkappa_Q$ and $\upsilon_Q$ are given, respectively, in \eqref{eq: K_*}, \eqref{eq: def Xi1}, \eqref{eq: def Xi2}, \eqref{eq:vcon} and \eqref{eq: upsilonQ}.
}
\rem{rem: univ correct}{
Note that the correction term of order $n^{-1/4}$ is not universal, as it depends on the potential $Q$. Indeed, although dimensionless, the factor $r_*\upsilon_Q/\gamma_Q$ cannot be factored out, since $\Xi_*^{(1)}$ and $\Xi_*^{(2)}$ are already independent of the potential.
}
\begin{proof}[Proof of \autoref{theo2 strong}]
In the proof of \autoref{theo2}, we saw in \eqref{eq: exp kernel 2} that the kernel admits the following expansion
\Eq{}{
 \frac{1}{\sqrt{\gamma_Q\,n}} K_{n}(z_1,z_2)=&\,\frac{1}{4}K_*(\xi_1,\xi_2)+\frac{1}{4 \varkappa_Q n^{1/4}}\Xi_*(\xi_1,\xi_2)+O(n^{-1/2}),
}
with
\Eq{}{
\Xi_*(\xi_1,\xi_2):=&\int_{-\infty}^{\infty}dv\, e^{i\Im(\xi_1-\xi_2)\frac{v^3}{12}} \frac{\exp\left(-\frac{v^4}{2}\left[ \mathbf{P}'\left(\frac{\Re(\xi_1)}{v}-1\right)+ \mathbf{P}'\left(\frac{\Re(\xi_2)}{v}-1\right)\right]\right)}{\mathcal{P}(v^4)}\abs{v}\\
  &\times\Bigg[\int_{-\infty}^{\infty}\left[\Phi_Q(vy+1,v)-vy\right]\frac{\exp[-v^4\mathbf{P}'\left(y\right)]}{\mathcal{P}(v^4)}dy-\frac{\Phi_Q(\xi_1,v)+\Phi_Q(\xi_2,v)}{2}\\
  &+\left(\frac{r_*\upsilon_Q}{\gamma_Q} + 4\right)\frac{v}{3}+i\Im(\xi_1-\xi_2)\left(\frac{r_*\upsilon_Q}{\gamma_Q} + 4\right)\frac{1}{2} \frac{v^4}{4!}\Bigg]
}
where $\Phi_Q$ is given in \eqref{eq: def phiQ}. Recalling the definition of $\mathbf{P}$ \eqref{eq: def bf P}, we have $L(y)=\frac{y^2}{2(3!)}+\mathbf{P}(y)$ and a simple expansion of the right-hand side of the above equation, using \eqref{eq: def L} and \eqref{eq: fct}, yields
\Eq{}{
\Xi_*(\xi_1,\xi_2)=\Xi_*^{(1)}(\xi_1,\xi_2)+\frac{r_*\upsilon_Q}{\gamma_Q} \Xi_*^{(2)}(\xi_1,\xi_2).
}
\end{proof}

\section*{Data availability}
This work does not have any associated data.
\section*{Acknowledgement}
We would like to thank the anonymous referees for their careful review of the manuscript and suggesting improvements. We are also grateful to Yacin Ameur, Sung-Soo Byun, Joakim Cronvall, Peter Forrester, Mario Kieburg, Arno Kuijlaars  and Aron Wennman for their valuable insights and comments. M.A. and S.L. acknowledge financial support from the International Research Training Group (IRTG) between KU Leuven and University of Melbourne and  Melbourne Research Scholarship of University of Melbourne. M.A. was also partially supported by the Australian Research Council Discovery Project DP250102552.
\section*{Conflict of interest}
The authors have no competing interests to declare that are relevant to the content of this article.
\bibliography{main}
\end{document}